\documentclass[12pt]{article}
\usepackage{latexsym}
\usepackage{amssymb}
\usepackage{amsmath}

\begin{document}


\newtheorem{theorem}{Theorem}[section]
\newtheorem{problem}{Problem}
\newtheorem{definition}{Definition}[section]
\newtheorem{lemma}{Lemma}[section]
\newtheorem{proposition}{Proposition}[section]
\newtheorem{corollary}{Corollary}[section]
\newtheorem{example}{Example}
\newtheorem{conjecture}{Conjecture}
\newtheorem{algorithm}{Algorithm}
\newtheorem{exercise}{Exercise}
\newtheorem{remarkk}{Remark}[section]

\newcommand{\be}{\begin{equation}}
\newcommand{\ee}{\end{equation}}
\newcommand{\bea}{\begin{eqnarray}}
\newcommand{\eea}{\end{eqnarray}}
\newcommand{\beq}[1]{\begin{equation}\label{#1}}
\newcommand{\eeq}{\end{equation}}
\newcommand{\beqn}[1]{\begin{eqnarray}\label{#1}}
\newcommand{\eeqn}{\end{eqnarray}}
\newcommand{\beaa}{\begin{eqnarray*}}
\newcommand{\eeaa}{\end{eqnarray*}}
\newcommand{\req}[1]{(\ref{#1})}

\newcommand{\lip}{\langle}
\newcommand{\rip}{\rangle}
\newcommand{\uu}{\underline}
\newcommand{\oo}{\overline}
\newcommand{\La}{\Lambda}
\newcommand{\la}{\lambda}
\newcommand{\eps}{\varepsilon}
\newcommand{\om}{\omega}
\newcommand{\Om}{\Omega}
\newcommand{\ga}{\gamma}
\newcommand{\rrr}{{\Bigr)}}
\newcommand{\qqq}{{\Bigl\|}}

\newcommand{\dint}{\displaystyle\int}
\newcommand{\dsum}{\displaystyle\sum}
\newcommand{\dfr}{\displaystyle\frac}
\newcommand{\bige}{\mbox{\Large\it e}}
\newcommand{\integers}{{\Bbb Z}}
\newcommand{\rationals}{{\Bbb Q}}
\newcommand{\reals}{{\rm I\!R}}
\newcommand{\realsd}{\reals^d}
\newcommand{\realsn}{\reals^n}
\newcommand{\NN}{{\rm I\!N}}
\newcommand{\DD}{{\rm I\!D}}
\newcommand{\LL}{{\rm I\!L}}
\newcommand{\degree}{{\scriptscriptstyle \circ }}
\newcommand{\dfn}{\stackrel{\triangle}{=}}
\def\complex{\mathop{\raise .45ex\hbox{${\bf\scriptstyle{|}}$}
     \kern -0.40em {\rm \textstyle{C}}}\nolimits}
\def\hilbert{\mathop{\raise .21ex\hbox{$\bigcirc$}}\kern -1.005em {\rm\textstyle{H}}} 
\newcommand{\RAISE}{{\:\raisebox{.6ex}{$\scriptstyle{>}$}\raisebox{-.3ex}
           {$\scriptstyle{\!\!\!\!\!<}\:$}}} 

\newcommand{\hh}{{\:\raisebox{1.8ex}{$\scriptstyle{\degree}$}\raisebox{.0ex}
           {$\textstyle{\!\!\!\! H}$}}}

\newcommand{\OO}{\won}
\newcommand{\calA}{{\cal A}}
\newcommand{\calB}{{\cal B}}
\newcommand{\calC}{{\cal C}}
\newcommand{\calD}{{\cal D}}
\newcommand{\calE}{{\cal E}}
\newcommand{\calF}{{\cal F}}
\newcommand{\calG}{{\cal G}}
\newcommand{\calH}{{\cal H}}
\newcommand{\calK}{{\cal K}}
\newcommand{\calL}{{\cal L}}
\newcommand{\calM}{{\cal M}}
\newcommand{\calO}{{\cal O}}
\newcommand{\calP}{{\cal P}}
\newcommand{\calX}{{\cal X}}
\newcommand{\calXX}{{\cal X\mbox{\raisebox{.3ex}{$\!\!\!\!\!-$}}}}
\newcommand{\calXXX}{{\cal X\!\!\!\!\!-}}
\newcommand{\gi}{{\raisebox{.0ex}{$\scriptscriptstyle{\cal X}$}
\raisebox{.1ex} {$\scriptstyle{\!\!\!\!-}\:$}}}
\newcommand{\intsim}{\int_0^1\!\!\!\!\!\!\!\!\!\sim}
\newcommand{\intsimt}{\int_0^t\!\!\!\!\!\!\!\!\!\sim}
\newcommand{\pp}{{\partial}}
\newcommand{\al}{{\alpha}}
\newcommand{\sB}{{\cal B}}
\newcommand{\sL}{{\cal L}}
\newcommand{\sF}{{\cal F}}
\newcommand{\sE}{{\cal E}}
\newcommand{\sX}{{\cal X}}
\newcommand{\R}{{\rm I\!R}}
\newcommand{\vp}{\varphi}
\newcommand{\N}{{\rm I\!N}}
\def\ooo{\lip}
\def\ccc{\rip}
\newcommand{\ot}{\hat\otimes}
\newcommand{\rP}{{\Bbb P}}
\newcommand{\bfcdot}{{\mbox{\boldmath$\cdot$}}}

\renewcommand{\varrho}{{\ell}}
\newcommand{\dett}{{\textstyle{\det_2}}}
\newcommand{\sign}{{\mbox{\rm sign}}}
\newcommand{\TE}{{\rm TE}}
\newcommand{\TA}{{\rm TA}}
\newcommand{\E}{{\rm E\,}}
\newcommand{\won}{{\mbox{\bf 1}}}
\newcommand{\Lebn}{{\rm Leb}_n}
\newcommand{\Prob}{{\rm Prob\,}}
\newcommand{\sinc}{{\rm sinc\,}}
\newcommand{\ctg}{{\rm ctg\,}}
\newcommand{\loc}{{\rm loc}}
\newcommand{\trace}{{\,\,\rm trace\,\,}}
\newcommand{\Dom}{{\rm Dom}}
\newcommand{\ifff}{\mbox{\ if and only if\ }}
\newcommand{\proof}{\noindent {\bf Proof:\ }}
\newcommand{\remark}{\noindent {\bf Remark:\ }}
\newcommand{\remarks}{\noindent {\bf Remarks:\ }}
\newcommand{\note}{\noindent {\bf Note:\ }}

\newcommand{\boldx}{{\bf x}}
\newcommand{\boldX}{{\bf X}}
\newcommand{\boldy}{{\bf y}}
\newcommand{\boldR}{{\bf R}}
\newcommand{\uux}{\uu{x}}
\newcommand{\uuY}{\uu{Y}}

\newcommand{\limn}{\lim_{n \rightarrow \infty}}
\newcommand{\limN}{\lim_{N \rightarrow \infty}}
\newcommand{\limr}{\lim_{r \rightarrow \infty}}
\newcommand{\limd}{\lim_{\delta \rightarrow \infty}}
\newcommand{\limM}{\lim_{M \rightarrow \infty}}
\newcommand{\limsupn}{\limsup_{n \rightarrow \infty}}

\newcommand{\ra}{ \rightarrow }

\newcommand{\ARROW}[1]
  {\begin{array}[t]{c}  \longrightarrow \\[-0.2cm] \textstyle{#1} \end{array} }

\newcommand{\AR}
 {\begin{array}[t]{c}
  \longrightarrow \\[-0.3cm]
  \scriptstyle {n\rightarrow \infty}
  \end{array}}

\newcommand{\pile}[2]
  {\left( \begin{array}{c}  {#1}\\[-0.2cm] {#2} \end{array} \right) }

\newcommand{\floor}[1]{\left\lfloor #1 \right\rfloor}

\newcommand{\mmbox}[1]{\mbox{\scriptsize{#1}}}

\newcommand{\ffrac}[2]
  {\left( \frac{#1}{#2} \right)}

\newcommand{\one}{\frac{1}{n}\:}
\newcommand{\half}{\frac{1}{2}\:}

\def\le{\leq}
\def\ge{\geq}
\def\lt{<}
\def\gt{>}

\def\squarebox#1{\hbox to #1{\hfill\vbox to #1{\vfill}}}
\newcommand{\qed}{\hspace*{\fill}
           \vbox{\hrule\hbox{\vrule\squarebox{.667em}\vrule}\hrule}\bigskip}

\title{Monge-Kantorovitch Measure  Transportation and Monge-Amp\`ere
  Equation  on Wiener Space}

\author{D. Feyel and A. S. \"Ust\"unel}
\date{ }
\maketitle
\noindent
{\bf Abstract:} 
{\small{Let $(W,\mu,H)$ be an abstract Wiener space assume two
    $\nu_i,\,i=1,2$  probabilities on $(W,\calB(W))$  
    {\footnote{ cf. Theorem \ref{monge-general} for the precise
        hypothesis about $\nu_1$ and $\nu_2$.}}. We give some conditions
    for the  Wasserstein distance between $\nu_1$ and $\nu_2$ 
    with respect to the  Cameron-Martin space  
    $$
    d_H(\nu_1,\nu_2)=\sqrt{\inf_\beta\int_{W\times W}|x-y|_H^2d\beta(x,y)} 
    $$ 
    to be finite, where the infimum is taken on  the set of probability
     measures $\beta$  on $W\times W$ whose first and second marginals
     are $\nu_1$ 
    and $\nu_2$. In this latter situation we prove the existence of  a
    unique (cyclically monotone) map 
    $T=I_W+\xi$, with $\xi:W\to H$, 
    such that $T$  maps $\nu_1$ to $\nu_2$. Besides, if $\nu_2\ll
    \mu${\footnote{In fact this hypothesis is too strong, cf. Theorem
       \ref{monge-general}.},  then   $T$ is stochastically 
    invertible, i.e., there exists $S:W\to W$ such that $S\circ T=I_W$
    $\nu_1$ a.s. and $T\circ S=I_W$ $\nu_2$ a.s. If $\nu_1=\mu$, then
    there exists  
    a $1$-convex function $\phi$ in the Gaussian Sobolev space
    $\DD_{2,1}$, such that $\xi=\nabla\phi$.   These
    results imply that the quasi-invariant transformations of the
    Wiener space with finite Wasserstein distance from $\mu$ can be
    written as the composition of a transport map $T$ and a rotation,
    i.e., a measure preserving map. We  give also
    1-convex sub-solutions and Ito-type  solutions  of the Monge-Amp\`ere
    equation on $W$.  
    }}   
\section{Introduction}
\label{intro}

In 1781, Gaspard Monge has published his celebrated memoire about the
most economical way of earth-moving \cite{Monge}. The configurations
of excavated earth and remblai were  modelized as two measures of
equal mass, say $\rho$ and $\nu$, that Monge had supposed absolutely
continuous with respect to the volume measure. Later Amp\`ere has
studied an analogous question about the electricity current in a media
with varying conductivity. In modern language of measure theory we can
express the problem in the following terms: let $W$ be a Polish space
on which are given two positive measures $\rho$ and $\nu$,  of finite,
equal mass. Let $c(x,y)$ be a cost function on $W\times W$, which is,
usually, assumed  positive. Does there  exist a map $T:W\to W$ such that
$T\rho=\nu$ and $T$ minimizes the integral 
$$
\int_W c(x,T(x))d\rho(x)
$$
between all such maps? The problem has been further studied by Appell
 \cite{App-1,App-2} and by Kantorovitch \cite{Kan}. Kantarovitch has
 succeeded to transform this highly nonlinear problem of Monge into a
 linear problem by replacing the search for $T$ with the search of a
 measure $\ga$ on $W\times W$ with marginals $\rho$ and $\nu$  such that the
 integral 
$$
\int_{W\times W}c(x,y)d\ga(x,y)
$$
is the minimum  of all the integrals 
$$
\int_{W\times W}c(x,y)d\beta(x,y)
$$
where $\beta$ runs in the set of measures on $W\times W$ whose
marginals are $\rho$ and $\nu$. Since then the problem adressed above
is called the  Monge problem and the quest of the optimal
measure is called the  Monge-Kantorovitch problem.

In this paper we study the Monge-Kantorovitch and the 
Monge   problem in the frame of an abstract  Wiener
space with a singular cost. In other words, let  $W$ be a separable
Fr\'echet space with its Borel sigma algebra 
$\calB(W)$ and assume that there is a separable  Hilbert space $H$ which is
injected densely and continuously into $W$, hence in general the
topology of $H$ is stronger than the topology induced by $W$. The cost
function $c:W\times W\to \reals_+\cup\{\infty\}$ is defined  as 
$$
c(x,y)=|x-y|_H^2\,,
$$
we suppose that $c(x,y)=\infty$ if $x-y$ does not belong to
$H$. Clearly, this
choice of the function $c$ is not arbitrary, in fact it is closely
related to Ito Calculus, hence also to the problems originating from 
Physics, quantum chemistry, large deviations, etc. Since
for all the interesting measures on $W$, the Cameron-Martin space is a
negligeable set, the cost function will be infinity very frequently. 
Let  $\Sigma(\rho,\nu)$ denote  the set  of probability measures 
on $W\times W$ with given marginals $\rho$ and $\nu$. It is a convex,  compact
 set under the weak topology $\sigma(\Sigma,C_b(W\times W))$. As
 explained above, the problem of Monge consists of finding a measurable map 
$T:W\to W$, called the optimal transport of $\rho$ to $\nu$, i.e.,
$T\rho=\nu$\footnote{We denote the
  push-forward of $\rho$ by $T$, i.e., the image of $\rho$ under $T$,  by $T\rho$.} which minimizes
the cost
$$
U\to  \int_W|x-U(x)|_H^2d\rho(x)\,,
$$
between all the maps $U:W\to W$ such that $U\rho=\nu$. The 
Monge-Kantorovitch  problem will   consist of
finding a  measure on $W\times W$,  which minimizes the function
$\theta\to J(\theta)$, defined by
\begin{equation}
\label{J-defn}
J(\theta)=\int_{W\times W}\big|x-y\big|_H^2d\theta(x,y)\,,
\end{equation}
where $\theta$ runs in $\Sigma(\rho,\nu)$. Note that 
$\inf\{J(\theta):\,\theta\in \Sigma(\rho,\nu)\}$ is the square of  Wasserstein
metric $d_H(\rho,\nu)$ with respect to the Cameron-Martin space $H$. 

Any solution $\ga$ of the
 Monge-Kantorovitch problem will give a solution to  the
Monge problem  provided that its support is included in
the graph of a map. Hence our work consists of realizing this program.
Although in  the finite dimensional case this problem is well-studied
in the path-breaking papers of Brenier \cite{BRE} and McCann \cite{Mc1,Mc2}
the things do not come up easily in our setting  and the  difficulty
is due to the 
fact that the cost function is not continuous with respect to the
Fr\'echet topology of $W$, for instance  the weak convergence of the
probability measures does not imply the convergence of the integrals  of
the cost function.  In other words the
function $|x-y|_H^2$ takes the value plus  infinity ``very
often''. On the other hand the results we obtain  seem to
have important applications to several problems of stochastic
analysis that we shall explain while enumerating the contents of the
paper. 

In Section \ref{preliminaries}, we explain some basic results about  the
functional analysis constructed on the Wiener space (cf., for instance
\cite{F-P,ASU}) and the probabilistic  theory of convex functions
recently developped in \cite{F-U1}. Section \ref{inequalities} is
devoted to the derivation of some inequalities which control the
Wasserstein distance. In particular, with the help of  the Girsanov
theorem, we give a very simple proof of an
inequality, initially  discovered by Talagrand (\cite{Tal}); this
facility gives already an idea about the efficiency of 
the infinite dimensional techniques for the Monge-Kantorovitch
problem{\footnote{In Section \ref{equation} we shall see another
    illustration of this phenomena.}}.  We indicate some simple
consequences of this inequality to 
control the measures of subsets of the Wiener space with respect to  second
moments of their  gauge functionals defined with  the
Cameron-Martin distance. These inequalities  are quite useful in the
theory of large deviations. Using a different representation of the
target measure, namely by constructing a flow of diffeomorphisms of
the Wiener space (cf. Chapter V of \cite{BOOK}) which maps the Wiener
measure to the target measure, we obtain also  a new control of the
Kantorovitch-Rubinstein 
metric of order one. The method we employ for this inequality
generalizes directly to a more general class of measures, namely those
for which one can define a reasonable divergence operator.

In Section \ref{gaussian}, we solve directly  the original problem of
Monge when the first measure is the Wiener measure and the second one is
given with a density, in such a way that the Wasserstein distance
between these two measures is finite. We prove the existence and the
uniqueness  of a
transformation of $W$ of the form $T=I_W+\nabla\phi$, where $\phi$ is a
$1$-convex function in the Gaussian Sobolev space $\DD_{2,1}$ such
that the measure $\ga=(I_W\times T)\mu$ is the unique solution of the
 problem of Monge-Kantorovitch.  This result gives
a new insight to the question  of representing an  integrable, positive
random variable whose expectation is unity,  as the Radon-Nikodym
derivative  of the image of the Wiener measure under a map which is a
perturbation of identity, a problem which has been studied by X. Fernique
and by one of us with M. Zakai (cf., \cite{Fer1,Fer2,BOOK}).  In \cite{BOOK},
Chapter II,  it is shown that such random 
variables are dense in $L^1_{1,+}(\mu)$ (the lower index $1$ means
that the expectations are equal to one), here we prove that this set
of random variables contains  the random variables who are at finite
Wasserstein distance from the Wiener measure. In fact even if this
distance is infinite, we show that there is a solution to this problem
if we  enlarge  $W$ slightly by taking $\NN\times W$. 

Section \ref{factorization} is devoted to the immediate
implications of the existence and the  uniqueness  of the solutions of 
Monge-Kantorovitch and Monge  problems constructed in Section \ref{gaussian}. 
Indeed the uniqueness  implies at once that the absolutely continuous
transformations of the Wiener space, at finite (Wasserstein) distance,
have a unique decomposition in the sense that they can be written as
the composition of a measure preserving map in the form of the
perturbation of identity with another one which is the perturbation of
identity with the Sobolev derivative of a $1$-convex function. This
means in particular that the class  of  $1$-convex functions is  as
basic  as the class  of adapted processes in the setting  of Wiener space.

In Section \ref{general-case} we prove the existence and the
uniqueness of solutions of the   Monge-Kantorovitch and Monge 
problems  for the measures which are at finite Wasserstein distance  from
each other. The fundamental  hypothesis  we use  is that the regular
conditional probabilities which are obtained by the disintegration of
one of the measures along the orthogonals of a sequence of regular, finite
dimensional projections vanish on the sets of co-dimension one. In
particular, this hypothesis  is satisfied if the measure under
question is absolutely continuous with respect to the Wiener measure.
The method  we use in this section is totally  different from the one
of Section \ref{gaussian}; it is  based on the notion of 
cyclic monotonicity of the supports of the regular conditional
probabilities obtained through some specific  disintegrations of the
optimal measures. The importance of cyclic monotonicity has first been
remarked by McCann and used abundently in \cite{Mc1} and in
\cite{G-Mc} for the finite dimensional case. Here the things are much
more complicated due to the singularity of the cost function, in
particular, contrary to the finite dimensional case,  the cyclic
monotonicity  is not compatible with the weak convergence of
probability measures.  A curious reader may ask why we did not treat
first the general case and then attack the subject of Section
\ref{gaussian}. The answer is twofold: even if we had done so, we
would have needed similar calculations as in Section \ref{gaussian} in
order to show the Sobolev regularity of the transport map, hence
concerning the volume, the order that we have chosen does not change
anything. Secondly, the construction used in  Section \ref{gaussian}
has  an interest by itself since it explains interesting  relations
between the transport map and its inverse and the optimal measure in a
more detectable situation, in this sense this construction  is rather
complementary to the material of Section \ref{general-case}.

Section \ref{equation} studies the Monge-Amp\`ere equation for the
measures which are absolutely continuous with respect to the Wiener
measure. First we briefly indicate the notion of  second order
Alexandroff derivative   and the Alexandroff version of the
Ornstein-Uhlenbeck operator applied to a $1$-convex function   in the
finite dimensional case. With the help of these observations, we write the
corresponding Jacobian using the modified Carleman-Fredholm
determinant which is natural in the infinite dimensional case (cf.,
\cite{BOOK}). Afterwards we attack the infinite dimensional case by
proving that the absolutely continuous part of the Ornstein-Uhlenbeck
operator applied to the finite rank conditional expectations of the
transport function is a submartingale which converges almost
surely. Hence the only difficulty lies in the calculation of the limit
of the Carleman-Fredholm determinants. Here we have a major difficulty
which originates from the pathology  of the
Radon-Nikodym derivatives of the vector measures with respect to a
scalar measure as explained in \cite{THOMAS}: in fact even if the
second order Sobolev  derivative of a Wiener function is a vector
measure with values in the space of Hilbert-Schmidt operators, its
absolutely continuous part has no reason to be Hilbert-Schmidt. Hence
the Carleman-Fredholm determinant may not exist, however 
due  to the $1$-convexity, the detereminants of the approximating
sequence are all with values in the interval $[0,1]$. Consequently  we
can construct the subsolutions with the help of the  Fatou lemma. 

Last but not the least, in section \ref{sub-equation}, we prove that
all these difficulties can be overcome thanks to the natural
renormalization of the Ito stochastic calculus. In fact using the Ito
representation theorem and the Wiener space analysis extended to the
distributions, cf. \cite{ASU-0},  we can give the explicit solution of the
Monge-Amp\`ere equation. This is a remarkable result in the sense that
such techniques do not exist in the finite dimensional case.


\section{Preliminaries and notations}
\label{preliminaries}
Let $W$ be a separable Fr\'echet space equipped with  a Gaussian
measure $\mu$ of zero mean whose support is the whole space. The
corresponding Cameron-Martin space is denoted by $H$. Recall that the
injection $H\hookrightarrow W$ is compact and its adjoint is the
natural injection $W^\star\hookrightarrow H^\star\subset
L^2(\mu)$. The triple $(W,\mu,H)$ is called 
an abstract Wiener space. Recall that $W=H$ if and only if $W$ is
finite dimensional. A subspace $F$ of $H$ is called regular if the
corresponding orthogonal projection 
has a continuous extension to $W$, denoted again  by the same letter.
It is well-known that there exists an increasing sequence of regular
subspaces $(F_n,n\geq 1)$, called total,  such that $\cup_nF_n$ is
dense in $H$ and in $W$. Let $\sigma(\pi_{F_n})${\footnote{For the notational
  simplicity, in the sequel we shall denote  it by  $\pi_{F_n}$.}}  be the
$\sigma$-algebra generated by $\pi_{F_n}$, then  for any  $f\in
L^p(\mu)$, the martingale  sequence 
$(E[f|\sigma(\pi_{F_n})],n\geq 1)$
converges to $f$ (strongly if 
$p<\infty$) in $L^p(\mu)$. Observe that the function
$f_n=E[f|\sigma(\pi_{F_n})]$ can be identified with a function on the
finite dimensional abstract Wiener space $(F_n,\mu_n,F_n)$, where
$\mu_n=\pi_n\mu$. 

Since the translations of $\mu$ with the elements of $H$ induce measures
equivalent to $\mu$, the G\^ateaux  derivative in $H$ direction of the
random variables is a closable operator on $L^p(\mu)$-spaces and  this
closure will be denoted by $\nabla$ cf.,  for example
\cite{F-P,ASU}. The corresponding Sobolev spaces 
(the equivalence classes) of the  real random variables 
will be denoted as $\DD_{p,k}$, where $k\in \NN$ is the order of
differentiability and $p>1$ is the order of integrability. If the
random variables are with values in some separable Hilbert space, say
$\Phi$, then we shall define similarly the corresponding Sobolev
spaces and they are denoted as $\DD_{p,k}(\Phi)$, $p>1,\,k\in
\NN$. Since $\nabla:\DD_{p,k}\to\DD_{p,k-1}(H)$ is a continuous and
linear operator its adjoint is a well-defined operator which we
represent by $\delta$. In the case of classical Wiener space, i.e.,
when $W=C(\reals_+,\reals^d)$, then $\delta$ coincides with the Ito
integral of the Lebesgue density of the adapted elements of
$\DD_{p,k}(H)$ (cf.\cite{ASU}). 

For any $t\geq 0$ and measurable $f:W\to \reals_+$, we note by
$$
P_tf(x)=\int_Wf\left(e^{-t}x+\sqrt{1-e^{-2t}}y\right)\mu(dy)\,,
$$
it is well-known that $(P_t,t\in \reals_+)$ is a hypercontractive
semigroup on $L^p(\mu),p>1$,  which is called the Ornstein-Uhlenbeck
semigroup (cf.\cite{F-P,ASU}). Its infinitesimal generator is denoted
by $-\calL$ and we call $\calL$ the Ornstein-Uhlenbeck operator
(sometimes called the number operator by the physicists). The
norms defined by 
\begin{equation}
\label{norm}
\|\phi\|_{p,k}=\|(I+\calL)^{k/2}\phi\|_{L^p(\mu)}
\end{equation}
are equivalent to the norms defined by the iterates of the  Sobolev
derivative $\nabla$. This observation permits us to identify the duals
of the space $\DD_{p,k}(\Phi);p>1,\,k\in\NN$ by $\DD_{q,-k}(\Phi')$,
with $q^{-1}=1-p^{-1}$, 
where the latter  space is defined by replacing $k$ in (\ref{norm}) by
$-k$, this gives us the distribution spaces on the Wiener space $W$
(in fact we can take as $k$ any real number). An easy calculation 
shows that, formally, $\delta\circ \nabla=\calL$, and this permits us
to extend the  divergence and the derivative  operators to the
distributions as linear,  continuous operators. In fact
$\delta:\DD_{q,k}(H\otimes \Phi)\to \DD_{q,k-1}(\Phi)$ and 
$\nabla:\DD_{q,k}(\Phi)\to\DD_{q,k-1}(H\otimes \Phi)$ continuously, for
any $q>1$ and $k\in \reals$, where $H\otimes \Phi$ denotes the
completed Hilbert-Schmidt tensor product (cf., for instance \cite{ASU}). 

Let us recall some facts from the convex analysis. Let $K$ be a 
Hilbert space, a subset $S$ of $K\times K$ is called cyclically
monotone if  any finite subset 
$\{(x_1,y_1),\ldots,(x_N,y_N)\}$ of
$S$ satisfies the following algebraic condition:
$$
\langle y_1,x_2-x_1\rangle+\langle y_2,x_3-x_2\rangle+\cdots+\langle
y_{N-1},x_N-x_{N-1}\rangle+\langle y_N,x_1-x_N\rangle\leq 0\,,
$$
where $\langle\cdot,\cdot\rangle$ denotes the inner product of
$K$. It turns out  that $S$ is
cyclically monotone if and only if 
$$
\sum_{i=1}^N(y_i,x_{\sigma(i)}-x_i)\leq 0\,,
$$
for any permutation $\sigma$ of $\{1,\ldots,N\}$ and for any finite
subset $\{(x_i,y_i):\,i=1,\ldots,N\}$ of $S$.
Note that  $S$ is
cyclically monotone if and only if any translate of it is cyclically
monotone.  By a theorem of Rockafellar,  any cyclically monotone set is
contained in the graph of the subdifferential  of a convex function in the
sense of convex analysis (\cite{ROC}) and even if the function may not
be unique its subdifferential  is unique. 

\noindent
Let now  $(W,\mu,H)$ be an abstract Wiener space;  a measurable  function
$f:W\to \reals\cup\{\infty\}$  is called $1$-convex if the map 
$$
h\to f(x+h)+\frac{1}{2}|h|_H^2=F(x,h)
$$
is convex on the Cameron-Martin space $H$ with values in
$L^0(\mu)$. Note that this notion is compatible with the
$\mu$-equivalence classes of random variables thanks to the
Cameron-Martin theorem. It is proven in \cite{F-U1} that 
this definition  is equivalent  the following condition:
  Let $(\pi_n,n\geq 1)$ be a sequence of regular, finite dimensional,
  orthogonal projections of 
  $H$,  increasing to the identity map
  $I_H$. Denote also  by $\pi_n$ its  continuous extension  to $W$ and
  define $\pi_n^\bot=I_W-\pi_n$. For $x\in W$, let $x_n=\pi_nx$ and
  $x_n^\bot=\pi_n^\bot x$.   Then $f$ is $1$-convex if and only if 
$$
x_n\to \frac{1}{2}|x_n|_H^2+f(x_n+x_n^\bot)
$$ 
is  $\pi_n^\bot\mu$-almost surely convex.


\section{Some Inequalities}
\label{inequalities}
\begin{definition}
Let $\xi$ and $\eta$ be two probabilities on $(W,\calB(W))$. We say
that a probability $\ga$ on $(W\times W,\calB(W\times W))$ is a
solution of the Monge-Kantorovitch problem associated to the
couple $(\xi,\eta)$ if the first marginal of $\ga$ is $\xi$, the
second one is $\eta$ and if 
$$
J(\ga)=\int_{W\times W}|x-y|_H^2d\ga(x,y)=\inf\left\{\int_{W\times
  W}|x-y|_H^2d\beta(x,y):\,\beta\in \Sigma(\xi,\eta)\right\}\,,
$$
where $\Sigma(\xi,\eta)$ denotes the set of all the probability
measures on $W\times W$ whose first and second marginals are
respectively $\xi$ and $\eta$. We shall denote the Wasserstein
distance between $\xi$ and $\eta$, which is  the positive
square-root of  this infimum, with $d_H(\xi,\eta)$. 
\end{definition}
\remark
Since the set of probability measures on $W\times W$ is weakly compact
and since the integrand in the definition is lower semi-continuous and
strictly convex,  the infimum in the definition is always  attained
even if the functional $J$ is identically infinity.

\noindent
The following result is an extension of an inequality due to Talagrand
\cite{Tal} and it gives a sufficient condition for the Wasserstein
distance to be finite:
\begin{theorem}
\label{ineq-thm}
Let $L\in \LL\log\LL(\mu)$ be a positive random variable with
$E[L]=1$ and let  $\nu$ be  the measure $d\nu=Ld\mu$.  We then  have  
\begin{equation}
\label{tal-ineq}
d_H^2(\nu,\mu)\leq 2E[L\log L]\,.
\end{equation}
\end{theorem}
\proof
Without loss of generality, we may suppose  that $W$ is equipped with
a filtration of sigma algebras in 
such a way that it becomes a classical Wiener space as
$W=C_0(\reals_+,\R^d)$. Assume first that  $L$ is  a strictly positive
and bounded  random variable. We can represent
it as 
$$
L=\exp\left[-\int_0^\infty(\dot{u}_s,dW_s)-\frac{1}{2}|u|_H^2\right]\,,
$$
where $u=\int_0^\cdot \dot{u}_sds$ is an $H$-valued, adapted random
variable. Define $\tau_n$ as
$$
\tau_n(x)=\inf\left\{t\in \reals_+:\,\int_0^t|\dot{u}_s(x)|^2ds>n\right\}\,.
$$
$\tau_n$ is a stopping time with respect to the canonical filtration
$(\calF_t,t\in \reals_+)$ of the Wiener process $(W_t,t\in \reals_+)$
and $\lim_n\tau_n=\infty$ almost surely. Define $u^n$ as
$$
u^n(t,x)=\int_0^t\won_{[0,\tau_n(x)]}(s)\dot{u}_s(x)ds\,.
$$
Let  $U_n:W\to W$ be  the map
$U_n(x)=x+u^n(x)$, then the Girsanov 
theorem says that $(t,x)\to U_n(x)(t)=x(t)+\int_0^t\dot{u}^n_sds$ is a
Wiener process under the measure $L_nd\mu$, where
$L_n=E[L|\calF_{\tau_n}]$. Therefore 
\beaa
E[L_n\log L_n]&=&E\left[L_n\,\left\{-\int_0^\infty
      (\dot{u}^n_s,dW_s)-\frac{1}{2}|u^n|_H^2\right\}\right]\\
&=&\frac{1}{2}E[L_n|u^n|_H^2]\\
&=&\frac{1}{2}E[L|u^n|_H^2]\,.
\eeaa
Define now the measure $\beta_n$ on $W\times W$ as 
$$
\int_{W\times W} f(x,y)d\beta_n(x,y)=\int_W f(U_n(x),x)L_n(x)d\mu(x)\,.
$$
Then the first marginal of $\beta_n$ is $\mu$ and the second one is
$L_n.\mu$. Consequently
\beaa
\lefteqn{\inf\left\{\int_{W\times
  W}|x-y|_H^2d\theta:\pi_1\theta=\mu,\,\pi_2\theta=L_n.\mu\right\}}\\
&\leq& \int_W|U_n(x)-x|_H^2L_nd\mu\\
&=&2E[L_n\log L_n]\,.
\eeaa
Hence we obtain 
$$
d_H^2(L_n.\mu,\mu)=J(\ga_n)\leq 2E[L_n\log L_n]\,,
$$
where $\ga_n$ is a solution of the Monge-Kantorovitch  problem in
$\Sigma(L_n.\mu,\mu)$. Let now $\ga$ be any cluster point of the
sequence $(\ga_n,n\geq 1)$, since $\ga\to J(\ga)$ is lower
semi-continuous with respect to the weak topology of probability
measures, we have 
\beaa
J(\ga)&\leq& \lim\inf_nJ(\ga_n)\\
&\leq& \sup_n2E[L_n\log L_n] \\
&\leq&2E[L\log L]\,,
\eeaa
since $\ga\in \Sigma(L.\mu,\mu)$, it follows that 
$$
d_H^2(L.\mu,\mu)\leq 2E[L\log L]\,.
$$
For the general case we stop the martingale $E[L|\calF_t]$
appropriately to obtain a bounded density $L_n$, then  replace it  by
$P_{1/n} L_n$ to improve the positivity, where $(P_t,t\geq 0)$ denotes
the Ornstein-Uhlenbeck semigroup. Then, from the Jensen inequality,
$$
E[P_{1/n}L_n\log P_{1/n}L_n]\leq E[L\log L]\,,
$$
therefore, using the same reasoning as above 
\beaa
d_H^2(L.\mu,\mu)&\leq&\lim\inf_nd_H^2(P_{1/n}L_n.\mu,\mu)\\
&\leq& 2E[L\log L]\,,
\eeaa
and this completes the proof.
\qed

\begin{corollary}
\label{tri-cor}
Assume that $\nu_i\,(i=1,2)$ have Radon-Nikodym  densities
$L_i\,(i=1,2)$ with respect to the Wiener measure $\mu$ which are in
$\LL\log\LL$. Then
$$
d_H(\nu_1,\nu_2)<\infty\,.
$$
\end{corollary}
\proof
This is a simple consequence of the triangle  inequality (cf. \cite{B-F}):
$$
d_H(\nu_1,\nu_2)\leq d_H(\nu_1,\mu)+d_H(\nu_2,\mu)\,.
$$
\qed

Let us give a simple application of the above result in the lines of
\cite{Mar}: 
\begin{corollary}
\label{iso-cor}
Assume that $A\in \calB(W)$ is any set of positive Wiener
measure. Define the $H$-gauge function of $A$ as 
$$
q_A(x)=\inf(|h|_H:\,h\in (A-x)\cap H)\,.
$$
Then we have 
$$
E[q_A^2]\leq 2\log\frac{1}{\mu(A)}\,,
$$
in other words
$$
\mu(A)\leq \exp\left\{-\frac{E[q_A^2]}{2}\right\}\,.
$$
Similarly if $A$ and $B$ are $H$-separated, i.e., if $A_\eps\cap
B=\emptyset$,  for some $\eps>0$, where  $A_\eps=\{x\in
W:\,q_A(x)\leq\eps\}$, then  
$$
\mu(A_\eps^c)\leq \frac{1}{\mu(A)}e^{-\eps^2/4}
$$
and consequently 
$$
\mu(A)\,\mu(B)\leq \exp\left(-\frac{\eps^2}{4}\right)\,.
$$
\end{corollary}
\remark We already know that, from the $0-1$--law,  $q_A$ is almost
surely finite, besides it satisfies $|q_A(x+h)-q_A(x)|\leq |h|_H$, hence
$E[\exp\la q_A^2]<\infty$ for any $\la<1/2$ (cf. \cite{BOOK}). In fact
all these assertions can also be proved with the technique used below.

\proof
Let $\nu_A$ be the measure defined by 
$$
d\nu_A=\frac{1}{\mu(A)}1_A d\mu\,.
$$
Let $\ga_A$ be the solution of the Monge-Kantorovitch problem, it is
easy to see that the support of $\ga_A$ is included in $W\times A$,
hence 
$$
|x-y|_H\geq \inf\{|x-z|_H:\,z\in A\}=q_A(x)\,,
$$
$\ga_A$-almost surely. This implies in particular that $q_A$ is almost
surely finite. It follows now from the inequality
(\ref{tal-ineq})
$$
E[q_A^2]\leq -2\log\mu(A)\,,
$$
hence the proof of the first inequality  follows. For the second let
$B=A_\eps^c$ and let $\ga_{AB}$ be the solution of the
Monge-Kantorovitch problem corresponding to $\nu_A,\nu_B$. Then we
have from the Corollary \ref{tri-cor}, 
$$
d^2_H(\nu_A,\nu_B)\leq -4\log\mu(A)\mu(B)\,.
$$
Besides the support of the measure $\ga_{AB}$ is in $A\times B$, hence
$\ga_{AB}$-almost surely $|x-y|_H\geq \eps$ and the proof follows.
\qed

\noindent
For the distance defined by 
$$
d_1(\nu,\mu)=\inf\left\{\int_{W\times
  W}|x-y|_H d\theta:\pi_1\theta=\mu,\,\pi_2\theta=\nu\right\}
$$
we have the following  control:

\begin{theorem}
\label{our-thm}
Let $L\in \LL_+^1(\mu)$ with $E[L]=1$. Then we have
\begin{equation}
\label{our-ineq}
d_1(L.\mu,\mu)\leq E\left[\left|(I+\calL)^{-1}\nabla L\right|_H\right]\,.
\end{equation}
\end{theorem}
\proof To prove the theorem we shall use a technique developed in
\cite{D-M}. 
Using the conditioning with respect to the sigma algebra
$V_n=\sigma\{\delta e_1,\ldots,\delta e_n\}$, where $(e_i,i\geq 1)$ is
a complete, orthonormal basis of $H$, 
we reduce the problem to the finite dimensional case. Moreover, we can
assume that $L$ is a smooth, strictly positive  function on
$\reals^n$. Define now $\sigma=(I+\calL)^{-1}\nabla L$ and 
$$
\sigma_t(x)=\frac{\sigma(x)}{t+(1-t)L}\,,
$$
for $t\in [0,1]$. Let $(\phi_{s,t}(x),s\leq t\in[0,1])$ be the flow of
diffeomorphisms defined by the following differential equation:
$$
\phi_{s,t}(x)=x-\int_s^t\sigma_\tau(\phi_{s,\tau}(x))d\tau\,.
$$
From the standart results (cf. \cite{BOOK}, Chapter V), it follows that
$x\to\phi_{s,t}(x)$ is Gaussian under the probability $\La_{s,t}.\mu$,
where 
$$
\La_{s,t}=\exp\int_s^t(\delta\sigma_\tau)(\phi_{s,\tau}(x))d\tau
$$
is the Radon-Nikodym density of $\phi_{s,t}^{-1}\mu$ with respect to $\mu$.
Define 
$$
H_s(t,x)=\La_{s,t}(x)\left\{t+(1-t)L\circ\phi_{s,t}(x)\right\}\,.
$$
It is easy to see that 
$$
\frac{d}{dt}H_s(t,x)=0
$$
for  $t\in(s,1)$. Hence the map $t\to H_s(t,x)$ is a constant, this
implies that 
$$
\La_{s,1}(x)=s+(1-s)L(x)\,.
$$
We have, as in the proof of Theorem \ref{ineq-thm}, 
\beaa
d_1(L.\mu,\mu)&\leq&E[|\phi_{0,1}(x)-x|_H\La_{0,1}]\\
&\leq&E\left[\La_{0,1}\int_0^1|\sigma_t(\phi_{0,t}(x))|_Hdt\right]\\
&=&E\left[\int_0^1\left|\sigma_t(\phi_{0,t}\circ\phi_{0,1}^{-1})(x)\right|_Hdt\right]\\
&=&E\left[\int_0^1\left|\sigma_t(\phi_{t,1}^{-1}(x))\right|_Hdt\right]\\
&=&E\left[\int_0^1|\sigma_t(x)|_H\La_{t,1}dt\right]\\
&=&E[|\sigma|_H]\,,
\eeaa
and the general case follows via the usual approximation procedure.
\qed


\section{Construction of the transport map}
\label{gaussian}
In this section we give  the construction of the transport map in  the
Gaussian case. We begin with the following lemma:
\begin{lemma}
\label{integ-lemma}
Let $(W,\mu,H)$ be an abstract Wiener space, assume that $f:W\to
\reals$ is a measurable function such that it is G\^ateaux  differentiable
in the direction of the Cameron-Martin space $H$, i.e., there exists
some $\nabla f:W\to H$ such that 
$$
f(x+h)=f(x)+\int_0^1(\nabla f(x+\tau h),h)_Hd\tau\,,
$$
$\mu$-almost surely, for any $h\in H$.
If $|\nabla f|_H\in L^2(\mu)$, then $f$ belongs to the Sobolev space
$\DD_{2,1}$. 
\end{lemma}
\proof 
Since $|\nabla |f||_H\leq |\nabla f|_H$, we can assume that $f$ is
positive. Moreover, for any $n\in \NN$, the function $f_n=\min(f,n)$
has also a G\^ateaux  derivative such that $|\nabla f_n|_H\leq |\nabla
f|_H$ $\mu$-almost surely. It follows from the Poincar\'e inequality
that the sequence $(f_n-E[f_n],n\geq 1)$ is bounded in $L^2(\mu)$,
hence it is also bounded in $L^0(\mu)$. Since $f$ is almost surely
finite, the sequence $(f_n, n\geq 1)$ is bounded in $L^0(\mu)$,
consequently the deterministic sequence $(E[f_n],n\geq 1)$ is also
bounded in $L^0(\mu)$. This means that $\sup_nE[f_n]<\infty$, 
hence  the monotone convergence theorem implies that  $E[f]<\infty$
and the proof is completed.
\qed

\begin{theorem}
\label{gaussian-case}
Let $\nu$ be the measure $d\nu=Ld\mu$, where $L$ is a positive random variable,
with $E[L]=1$. Assume that $d_H(\mu,\nu)<\infty$ (for instance 
$L\in \LL\log \LL$). Then there exists a  $1$-convex function $\phi\in
\DD_{2,1}$, unique upto a constant,  such that  the map
$T=I_W+\nabla \phi$ is the unique solution of the original problem of
Monge. Moreover, its graph supports  the  unique
solution of the 
Monge-Kantorovitch problem $\ga$. Consequently    
$$
(I_W\times T)\mu=\ga
$$
In particular  $T$ maps $\mu$ to $\nu$ and  $T$ is almost surely
invertible, i.e., there exists some $T^{-1}$ such that $T^{-1}\nu=\mu$
and that 
\beaa
1&=&\mu\left\{x:\,T^{-1}\circ T(x)=x\right\}\\
&=&\nu\left\{y\in W:\,T\circ T^{-1}(y)=y\right\}\,.
\eeaa
\end{theorem}
\proof Let $(\pi_n,n\geq 1)$ be a sequence of regular, finite
dimensional orthogonal projections of $H$ increasing to $I_H$. Denote
their continuous extensions to $W$ by the same letters. For $x\in W$,
we define $\pi_n^\bot x=:x_n^\bot=x-\pi_nx$. Let  $\nu_n$ be  the measure
$\pi_n\nu$. Since $\nu$ is absolutely continuous with respect to
$\mu$, $\nu_n$ is absolutely continuous with respect to
$\mu_n:=\pi_n\mu$ and 
$$
\frac{d\nu_n}{d\mu_n}\circ\pi_n=E[L|V_n]=:L_n\,,
$$
where $V_n$ is the sigma algebra  $\sigma(\pi_n)$ and the conditional
expectation is taken with 
respect to $\mu$. On the space $H_n$, the  Monge-Kantorovitch
problem, which consists of finding the probability measure which
realizes the following infimum 
$$
d_H^2(\mu_n,\nu_n)=\inf\left\{J(\beta):\,\beta\in M_1(H_n\times
H_n)\,,p_1\beta=\mu_n,p_2\beta=\nu_n\right\} 
$$
where
$$
J(\beta)=\int_{H_n\times H_n}|x-y|^2d\beta(x,y)\,,
$$
has a unique solution $\ga_n$, where $p_i,\,i=1,2$ denote  the
projections  $(x_1,x_2)\to x_i,\,i=1,2$ from $H_n\times H_n$ to $H_n$
and $M_1(H_n\times H_n)$ denotes the set of probability measures on
$H_n\times H_n$.  The measure $\ga_n$ may be regarded
as a measure on $W\times W$, by taking its image under the injection
$H_n\times H_n\hookrightarrow W\times W$ which  we shall denote again by
$\ga_n$. It results from the finite dimensional results of Brenier and
of McCann(\cite{BRE}, \cite{Mc1}) that 
there are two convex continuous  functions (hence almost everywhere
differentiable) $\Phi_n$ and $\Psi_n$ on 
$H_n$ such that 
$$
\Phi_n(x)+\Psi_n(y)\geq (x,y)_H
$$
for all $x,y\in H_n$ and that 
$$
\Phi_n(x)+\Psi_n(y)= (x,y)_H
$$
$\ga_n$-almost everywhere. Hence the 
support of $\ga_n$ is included in the graph of the derivative
$\nabla\Phi_n$  of $\Phi_n$, hence $\nabla\Phi_n\mu_n=\nu_n$ and the inverse of
$\nabla\Phi_n$ is equal to $\nabla\Psi_n$. Let 
\beaa
\phi_n(x)&=&\Phi_n(x)-\frac{1}{2}|x|_H^2\\
\psi_n(y)&=&\Psi_n(y)-\frac{1}{2}|y|_H^2\,.
\eeaa
Then $\phi_n$ and $\psi_n$ are $1$-convex functions and they satisfy
the following relations:
\begin{equation}
\label{+-ineq}
\phi_n(x)+\psi_n(y)+\frac{1}{2}|x-y|_H^2\geq 0\,,
\end{equation}
for all $x,y\in H_n$ and 
\begin{equation}
\label{=-eq}
\phi_n(x)+\psi_n(y)+\frac{1}{2}|x-y|^2_H=0\,,
\end{equation}
$\ga_n$-almost everywhere. From what we have said above, it follows
that $\ga_n$-almost surely  $y=x+\nabla\phi_n(x)$, consequently
\begin{equation}
\label{energy-iden}
J(\ga_n)=E[|\nabla\phi_n|_H^2]\,.
\end{equation}
Let $q_n:W\times W\to H_n\times H_n$ be defined as
$q_n(x,y)=(\pi_nx,\pi_ny)$. If $\ga$ is any solution of the 
Monge-Kantorovitch problem, then $q_n\ga\in \Sigma(\mu_n,\nu_n)$,
hence
\begin{equation}
\label{1st-estimate}
J(\ga_n)\leq J(q_n\ga)\leq J(\ga)=d_H^2(\mu,\nu)\,.
\end{equation}
Combining the relation (\ref{energy-iden}) with the inequality
(\ref{1st-estimate}), we obtain the following bound
\bea
\label{control-ineq}
\sup_nJ(\ga_n)&=&\sup_nd_H^2(\mu_n,\nu_n)\nonumber\\
&=&\sup_nE[|\nabla\phi_n|_H^2]\nonumber\\
&\leq&d_H^2(\mu,\nu)=J(\ga)\,.
\eea
For $m\leq n$, $q_m\ga_n\in\Sigma(\mu_m,\nu_m)$, hence  we should have 
\beaa
J(\ga_m)&=&\int_{W\times W}|\pi_mx-\pi_my|_H^2d\ga_m(x,y)\\
&\leq&\int_{W\times W}|\pi_mx-\pi_my|_H^2d\ga_n(x,y)\\
&\leq&\int_{W\times W}|\pi_nx-\pi_ny|_H^2d\ga_n(x,y)\\
&=&\int_{W\times W}|x-y|_H^2d\ga_n(x,y)\\
&=&J(\ga_n)\,,
\eeaa
where the third equality follows from the fact that we have denoted
the $\ga_n$ on $H_n\times H_n$ and its image in $W\times W$ by the
same letter.
Let now $\ga$ be a weak cluster point of the sequence of measures
$(\ga_n,n\geq 1)$, where the word `` weak''\footnote{To prevent the
  reader against the trivial errors let us emphasize  that $\ga_n$ is
  not the projection of $\ga$ on $W_n\times W_n$.} refers to the weak
convergence of 
measures on $W\times W$. Since $(x,y)\to |x-y|_H$ is lower
semi-continuous, we have 
\beaa
J(\ga)&=&\int_{W\times W}|x-y|_H^2d\ga(x,y)\\
&\leq&\lim\inf_n\int_{W\times W}|x-y|_H^2d\ga_n(x,y)\\
&=&\lim\inf_nJ(\ga_n)\\
&\leq&\sup_nJ(\ga_n)\\
&\leq& J(\ga)=d_H^2(\mu,\nu)\,,
\eeaa
from the relation (\ref{control-ineq}). Consequently 
\be
\label{lim-eq}
J(\ga)=\lim_nJ(\ga_n)\,.
\ee
Again from (\ref{control-ineq}), if we replace $\phi_n$ with $\phi_n-E[\phi_n]$
and $\psi_n$ with $\psi_n+E[\phi_n]$ we obtain a bounded sequence
$(\phi_n,n\geq 1)$ in $\DD_{2,1}$, in particular it is bounded in the
space $L^2(\ga)$ if we inject it into  latter  by $\phi_n(x)\to
\phi_n(x)\otimes 1(y)$. Consider now the sequence of the 
positive, lower semi-continuous functions $(F_n,n\geq 1)$ defined on
$W\times W$ as   
$$
F_n(x,y)=\phi_n(x)+\psi_n(y)+\frac{1}{2}|x-y|_H^2\,.
$$
We have, from the relation (\ref{=-eq})
\beaa
\int_{W\times
  W}F_n(x,y)d\ga(x,y)&=&\int_W\phi_nd\mu+\int_W\psi_n(y)d\nu+\frac{1}{2}J(\ga)
\\ 
&=&\frac{1}{2}\left(J(\ga)-J(\ga_n)\right)\to 0\,.
\eeaa
Consequently the sequence $(F_n,n\geq 1)$ converges to zero in
$L^1(\ga)$, therefore it is uniformly integrable. Since $(\phi_n,n\geq
1)$ is uniformly integrable as explained above and since $|x-y|^2$ has
a finite expectation with respect to $\ga$, it follows that
$(\psi_n,n\geq 1)$ is also uniformly integrable in $L^1(\ga)$ hence
also in $L^1(\nu)$. Let $\phi'$ be a weak cluster point of
$(\phi_n,n\geq 1)$, then there exists a sequence $(\phi_n',n\geq 1)$
whose elements are the convex combinations of some  elements of
$(\phi_k,k\geq n)$ such that $(\phi_n',n\geq 1)$ converges in the norm
topology of $\DD_{2,1}$ and $\mu$-almost everywhere. Therefore the
 sequence $(\psi_n',n\geq 1)$, constructed from $(\psi_k,k\geq n)$,
 converges in $L^1(\nu)$ and $\nu$-almost surely. Define $\phi$ and
 $\psi$ as  
\beaa
\phi(x)&=&\lim\sup_n\phi_n'(x)\\
\psi(y)&=&\lim\sup_n\psi_n'(y)\,,
\eeaa
hence we have 
$$
G(x,y)=\phi(x)+\psi(y)+\frac{1}{2}|x-y|_H^2\geq 0
$$
for all $(x,y)\in W\times W$, also the  equality holds $\ga$-almost
everywhere. Let now $h$ be any element of $H$, since $x-y$ is in $H$
for $\ga$-almost all $(x,y)\in W\times W$, we have 
$$
|x+h-y|_H^2=|x-y|_H^2+|h|_H^2+2(h,x-y)_H
$$
$\ga$-almost surely. Consequently 
$$
\phi(x+h)-\phi(x)\geq -(h,x-y)_H-\frac{1}{2}|h|_H^2
$$
$\ga$-almost surely and  this implies that 
$$
y=x+\nabla \phi(x)
$$
$\ga$-almost everywhere. Define now the map $T:W\to W$ as
$T(x)=x+\nabla\phi(x)$, then 
\beaa
\int_{W\times W}f(x,y)d\ga(x,y)&=&\int_{W\times W}f(x,T(x))d\ga(x,y)\\
&=&\int_{W}f(x,T(x))d\mu(x)\,,
\eeaa
for any $f\in C_b(W\times W)$, consequently $(I_W\times T)\mu=\ga$, in
particular $T\mu=\nu$.

Let us notice
that any weak cluster point of $(\phi_n,n\geq 1)$, say $\tilde{\phi}$,
satisfies 
$$
\nabla\tilde{\phi}(x)=y-x
$$
$\ga$-almost surely, hence $\mu$-almost surely we have
$\tilde{\phi}=\phi$. This implies that $(\phi_n,n\geq 1)$ has a unique
cluster point $\phi$, consequently the sequence $(\phi_n,n\geq 1)$
converges weakly in $\DD_{2,1}$ to $\phi$. Besides we have
\beaa
\lim_n\int_W|\nabla\phi_n|_H^2d\mu&=&\lim_nJ(\ga_n)\\
&=&J(\ga)\\
&=&\int_{W\times W}|x-y|_H^2d\ga(x,y)\\
&=&\int_W|\nabla\phi|_H^2d\mu\,,
\eeaa
hence $(\phi_n,n\geq 1)$ converges to $\phi$ in the norm topology of
$\DD_{2,1}$.  Let us recapitulate what we
have done till here: we have taken an arbitrary optimal $\ga\in
\Sigma(\mu,\nu)$ and an arbitrary cluster point $\phi$ of
$(\phi_n,n\geq 1)$ and we have proved that $\ga$ is carried by the
graph of $T=I_W+\nabla\phi$. This implies that $\ga$ and $\phi$ are unique
and  that  the sequence $(\ga_n,n\geq 1)$ has a unique cluster point $\ga$.

Certainly  $(\psi_n,\geq 1)$
converges also in the norm topology of $L^1(\nu)$.
Moreover, from the finite dimensional situation, we have
$\nabla\phi_n(x)+\nabla\psi_n(y)=0$ $\ga_n$-almost everywhere. Hence
$$
E_\nu[|\nabla\psi_n|_H^2]=E[|\nabla\phi_n|_H^2]
$$
this implies the
boundedness of $(\nabla\psi_n,n\geq 1)$ in $L^2(\nu,H)$ (i.e.,
$H$-valued functions). 
To complete the proof
we have to show that, for some measurable, $H$-valued map, say
$\eta$, it holds that 
$x=y+\eta(y)$ $\ga$-almost surely. For this let 
$F$ be a finite dimensional, regular subspace of $H$ and denote by
$\pi_F$ the projection operator onto $F$ which is continuously
extended to $W$, put $\pi_F^\bot=I_W-\pi_F$. We have
$W=F\oplus F^\bot$, with $F^\bot=\ker \pi_F=\pi_F^\bot(W)$. Define the
measures $\nu_F=\pi_F(\nu)$ and $\nu_F^\bot=\pi_F^\bot(\nu)$. From the
construction of $\psi$, we know that, for 
any $v\in F^\bot$, the partial map $u\to \psi(u+v)$ is $1$-convex on
$F$. Let also $A=\{y\in W:\,\psi(y)<\infty\}$, then $A$ is a Borel set
with $\nu(A)=1$ and it is easy to see that, for  $\nu_F^\bot$-almost all
$v\in F^\bot$, one has 
$$
\nu(A|\pi_F^\bot=v)>0\,.
$$
It then follows from Lemma 3.4 of \cite{F-U1}, and from the fact that
the regular conditional probability $\nu(\cdot\,|\pi_F^\bot=v)$ is
absolutely continuous with respect to the Lebesgue measure of $F$,
that $u\to \psi(u+v)$ is $\nu(\cdot\,|\pi_F^\bot=v)$-almost everywhere
differentiable on $F$ for $\nu_F^\bot$-almost
all $v\in F^\bot$. It then follows that, $\nu$-almost surely,  $\psi$
is differentiable in the directions of $F$, i.e., there exists
$\nabla_F\psi\in F$ $\nu$-almost surely. Since we also have 
$$
\psi(y+k)-\psi(y)\geq (x-y,k)_H-\frac{1}{2}|k|_H^2\,,
$$
we obtain, $\ga$-almost surely
$$
(\nabla_F\psi(y),k)_H=(x-y,k)_H\,,
$$
for any $k\in F$. Consequently 
$$
\nabla_F\psi(y)=\pi_F(x-y)
$$
$\ga$-almost surely. Let now $(F_n,n\geq 1)$ be a total, increasing
sequence of regular subspaces of $H$, we have a sequence
$(\nabla_n\psi,n\geq 1)$ bounded in $L^2(\nu)$ hence also bounded in
$L^2(\ga)$. Besides $\nabla_n\psi(y)=\pi_nx-\pi_ny$ $\ga$-almost
surely. Since $(\pi_n(x-y),n\geq 1)$ converges in $L^2(\ga,H)$,
$(\nabla_n\psi,n\geq 1)$ converges in the norm topology of
$L^2(\ga,H)$. Let us denote  this limit by  $\eta$, 
then  we have $x=y+\eta(y)$ $\ga$-almost surely. Note that, since
$\pi_n\eta=\nabla_n\psi$, we can even write in a weak sense that
$\eta=\nabla\psi$.  If we define $T^{-1}(y)=y+\eta(y)$, we
see that 
\beaa
1&=&\ga\{(x,y)\in W\times W:T\circ T^{-1}(y)=y\}\\
&=&\ga\{(x,y)\in W\times W:T^{-1}\circ T(x)=x\}\,,  
\eeaa
and this completes the proof of the theorem.
\qed

\begin{remarkk}
\label{nu-closed}
 Assume that the operator $\nabla$ is closable with respect to $\nu$,
then we have $\eta=\nabla\psi$. In particular, if  $\nu$ and $\mu$ are
equivalent, then we have   
$$
T^{-1}=I_W+\nabla\psi\,,
$$
where is $\psi$ is  a $1$-convex function.
\end{remarkk}
\begin{remarkk}
Assume that $L\in \LL_+^1(\mu)$, with $E[L]=1$ and let $(D_k,k\in
\NN)$ be a measurable partition of $W$ such that on each $D_k$, $L$ is
bounded. Define $d\nu=L\,d\mu$ and $\nu_k=\nu(\cdot|D_k)$. It follows
 from Theorem \ref{ineq-thm}, that $d_H(\mu,\nu_k)<\infty$. 
Let then  $T_k$ be the map constructed in Theorem \ref{gaussian-case}
satisfying $T_k\mu=\nu_k$. Define  $n(dk)$ as  the probability
distribution on $\NN$ given by  $n\left(\{k\}\right)=\nu(D_k),\,k\in
\NN$. Then we have  
$$
\int_W f(y)d\nu(y)=\int_{W\times \NN}f(T_k(x))\mu(dx)n(dk)\,.
$$
A  similar result is given in  \cite{Fer2}, the difference with that of
above lies in the fact that we have a more precise information about
the probability space on which $T$ is defined.
\end{remarkk}


\section{Polar factorization of the absolutely continuous
  transformations of the Wiener space}
\label{factorization}
Assume that $V=I_W+v:W\to W$ be an absolutely continuous  transformation and
let $L\in \LL_+^1(\mu)$ be the Radon-Nikodym derivative of $V\mu$ with
respect to $\mu$. Let
$T=I_W+\nabla\phi$ be the 
transport map such that $T\mu=L.\mu$. Then it is easy to see that the
map $s=T^{-1}\circ V$ is a rotation, i.e., $s\mu=\mu$
(cf. \cite{BOOK})  and it can be represented as $s=I_W+\alpha$. In
particular we have  
\begin{equation}
\label{iden-1}
\alpha+\nabla\phi\circ s=v\,.
\end{equation}
Since $\phi$ is a $1$-convex map, we have $h\to
\frac{1}{2}|h|_H^2+\phi(x+h)$ is almost surely convex (cf.\cite{F-U1}). Let
$s'=I_W+\alpha'$ be another rotation with $\alpha':W\to H$. By the
$1$-convexity of $\phi$, we have 
$$
\frac{1}{2}|\alpha'|_H^2+\phi\circ s'\geq
\frac{1}{2}|\alpha|_H^2+\phi\circ s
+\Bigl(\alpha+\nabla\phi\circ s,\alpha'-\alpha\Bigr)_H\,,
$$
$\mu$-almost surely. 
Taking the expectation of both sides,  using the fact that $s$ and
$s'$ preserve the Wiener measure $\mu$ and the identity (\ref{iden-1}),
we obtain  
$$
E\left[\half|\alpha|_H^2-(v,\alpha)_H\right]\leq
E\left[\frac{1}{2}|\alpha'|_H^2-(v,\alpha')_H\right] \,.
$$
Hence we have proven the existence part of the following 
\begin{proposition}
Let $\mathcal R_2$ denote the subset of $L^2(\mu,H)$ whose elements
are defined by the property that  $x\to x+\eta(x)$ is a rotation,
i.e., it preserves the Wiener measure. Then $\alpha$ is the unique
element of $\mathcal R_2$ which minimizes the functional
$$
\eta\to M_v(\eta)=E\left[\frac{1}{2}|\eta|_H^2-(v,\eta)_H\right]\,.
$$
\end{proposition}
\proof
To show the uniqueness, assume that $\eta\in \mathcal R_2$ be another
map minimizing $J_v$. Let $\beta$ be the measure on $W\times W$,
defined as 
$$
\int_{W\times W} f(x,y)d\beta(x,y)=\int_W f(x+\eta(x),V(x))d\mu\,.
$$
Then the first marginal of $\beta$ is $\mu$ and the second marginal is
$L.\mu$. Since $\ga=(I_W\times T)\mu$ is the unique solution of the
Monge-Kantorovitch problem, we should have 
$$
\int|x-y|_H^2 d\beta(x,y)>\int
|x-y|_H^2d\ga(x,y)=E[|\nabla\phi|_H^2]\,.
$$
However we have 
\beaa
\int_{W\times W}|x-y|_H^2 d\beta(x,y)&=&E\left[|v-\eta|_H^2\right]\\
&=&E\left[|v|_H^2\right]+2M_v(\eta)\\
&=&E\left[|v|_H^2\right]+2M_v(\alpha)\\
&=&E\left[|v-\alpha|_H^2\right]\\
&=&E\left[|\nabla \phi\circ s|_H^2\right]\\
&=&E\left[|\nabla \phi|_H^2\right]\\
&=&\int_{W\times W}|x-y|_H^2 d\ga(x,y)\\
&=&J(\ga)
\eeaa
and this  gives a contradiction to the uniqueness of $\ga$.
\qed

The following theorem, whose proof is rather easy, gives a better
understanding of  the structure of  absolutely continuous
transformations of the Wiener measure:
\begin{theorem}
\label{pol-fact}
Assume that $U:W\to W$ be a measurable  map and $L\in \LL\log \LL$ a
positive random variable with $E[L]=1$. Assume that the measure
$\nu=L\cdot\mu$ is a Girsanov measure for $U$, i.e., that one has 
$$
E[f\circ U\,L]=E[f]\,,
$$
for any $f\in C_b(W)$. Then there exists a unique map
$T=I_W+\nabla\phi$ with $\phi\in \DD_{2,1}$ is $1$-convex, and a measure
preserving transformation $R:W\to W$ such that $U\circ T=R$
$\mu$-almost surely and $U=R\circ T^{-1}$ $\nu$-almost surely.
\end{theorem}
\proof
By Theorem \ref{gaussian-case} there is a unique map
$T=I_W+\nabla\phi$, with $\phi\in\DD_{2,1}$,  $1$-convex such that $T$
transports $\mu$ to $\nu$. Since $U\nu=\mu$, we have 
\beaa
E[f\circ U\,L]&=&E[f\circ U\circ T]\\
&=&E[f]\,.
\eeaa
Therefore $x\to U\circ T(x)$ preserves the measure $\mu$. The rest is
obvious since $T^{-1}$ exists $\nu$-almost surely.
\qed

Another version of Theorem \ref{pol-fact} can be announced as follows:
\begin{theorem}
\label{pol-fact2}
Assume that $Z:W\to W$ is a measurable map such that $Z\mu\ll\mu$,
with $d_H(Z\mu,\mu)<\infty$. Then $Z$ can be decomposed as 
$$
Z=T\circ s\,,
$$
where $T$ is the unique  transport map of the Monge-Kantorovitch
problem for $\Sigma(\mu,Z\mu)$ and $s$ is a rotation.
\end{theorem}
\proof
Let $L$ be the Radon-Nikodym derivative of $Z\mu$ with respect to 
$\mu$. We have, from Theorem \ref{gaussian-case},  
\beaa
E[f]&=&E[f\circ T^{-1}\circ T]\\
&=&E[f\circ T^{-1}\,L]\\
&=&E[f\circ T^{-1}\circ Z]\,,
\eeaa
for any $f\in C_b(W)$. Hence $T^{-1}\circ Z=s$ is a rotation. Since
$T$ is uniquely defined, $s$ is also uniquely defined.
\qed

Although the following result is a translation of the results of this section,
 it is interesting from the point of view of  stochastic differential
 equations:
\begin{theorem}
\label{weak-sde}
Let  $(W,\mu,H)$ be  the standard  Wiener space on $\reals^d$,
i.e., $W=C(\reals_+,\reals^d)$. Assume that there exists a
probability $P\ll\mu$ which is the weak solution of the stochastic
differential equation  
$$
dy_t=dW_t+b(t,y)dt\,,
$$
such that $d_H(P,\mu)<\infty$. Then there exists a process $(T_t,t\in
\reals_+)$ which is  a pathwise  solution of some stochastic
differential equation whose law is equal to $P$.
\end{theorem}
\proof
Let $T$ be the transport map constructed in Theorem
\ref{gaussian-case} corresponding to $dP/d\mu$. Then it has an inverse
$T^{-1}$ such that $\mu\{T^{-1}\circ T(x)=x\}=1$. Let $\phi$ be the
$1$-convex function  such that $T=I_W+\nabla \phi$ and denote by
$(D_s\phi,s\in \reals_+)$ the representation of $\nabla\phi$ in
$L^2(\reals_+,ds)$. Define $T_t(x)$ as the trajectory $T(x)$ evaluated
at $t\in \reals_+$. Then it is easy to see that $(T_t,t\in \reals_+)$
satifies the stochastic differential equation
$$
T_t(x)=W_t(x)+\int_0^t l(s,T(x))ds\,\,,\,\,t\in \reals_+\,,
$$
where $W_t(x)=x(t)$ and $l(s,x)=D_s\phi\circ T^{-1}(x)$.
\qed


\section{Construction and uniqueness  of the transport map in the
  general case} 
\label{general-case}
In this section we call optimal every probability measure{\footnote{In
  fact the results of this section are essentially true for the
  bounded, positive measures.}} $\ga$ on
$W\times W$ such that $J(\ga)<\infty$ and that $J(\ga)\leq J(\theta)$
for every other probability $\theta$ having the same marginals as
those of $\ga$. We recall that a  finite dimensional subspace $F$  of $W$ is 
called regular if the corresponding  projection is
continuous. Similarly a finite dimensional projection of $H$ is called
regular if it has a continuous extension to $W$.

We begin with the following lemma which answers all kind of questions
of measurability that we may encounter in the sequel:
\begin{lemma}
\label{meas-lemma}
Consider two uncountable Polish spaces $X$ and $T$. Let $t\to \ga_t$
be a Borel family of probabilities on $X$ and let $\calF$ be a
separable sub-$\sigma$-algebra of the Borel $\sigma$-algebra $\calB$
of $X$. Then there exists a Borel kernel 
$$
N_tf(x)=\int_Xf(y)N_t(x,dy)\,,
$$
such that, for any bounded Borel function $f$ on $X$, the following
properties hold true:
\begin{enumerate}
\item[i)]$(t,x)\to N_tf(x)$ is Borel measurable on $T\times X$.
\item[ii)]For any $t\in T$, $N_tf$ is an $\calF$-measurable version of
  the conditional expectation $E_{\ga_t}[f|\calF]$.
\end{enumerate}
\end{lemma}
\proof
Assume first that $\calF$ is finite, hence it is generated by a finite
partition $\{A_1,\ldots,A_k\}$. In this case it suffices to take 
$$
N_tf(x)=\sum_{i=1}^k\frac{1}{\ga_t(A_i)}\left(\int_{A_i}fd\ga_t\right)\,1_{A_i}(x)\,\,\left({\mbox{
    with }}0=\frac{0}{0}\right)\,.
$$
For the general case, take an increasing sequence $(\calF_n,n\geq 1)$
of finite sub-$\sigma$-algebras whose union generates $\calF$. Without
loss of generality we can assume that $(X,\calB)$ is the Cantor set
(Kuratowski Theorem, cf., \cite{Del-M}). Then for every clopen set
(i.e., a set which is closed and open at the same time) $G$ and any
$t\in T$, the sequence $(N_t^n1_G,n\geq 1)$ converges $\ga_t$-almost
everywhere. Define
$$
H_G(t,x)=\limsup_{m,n\to\infty}|N^n_t1_G(x)-N^m_t1_G(x)|\,.
$$
$H_G$ is a Borel function on $T\times X$ which vanishes $\ga_t$-almost
all $x\in X$, besides, for any $t\in T$,  $x\to H_G(t,x)$ is
$\calF$-measurable. As there exist only countably many clopen sets in
$X$, the function 
$$
H(t,x)=\sup_GH_G(t,x)
$$
inherits all the  measurability properties. Let $\theta$ be any
probability on $X$, for any clopen $G$, define
$$
\begin{array}{cllrl}
N_t1_G(x)&=&\lim_nN_t^n1_G(x)&{\mbox{ if }}& H(t,x)=0\,,\\
&=&\theta(G)&{\mbox{ if }}& H(t,x)>0\,.
\end{array}
$$
Hence, for any $t\in T$, we get an additive measure on the Boolean
algebra of clopen sets of $X$. Since such a measure is
$\sigma$-additive and extends uniquely as a $\sigma$-additive measure
on $\calB$, the proof is completed.
\qed

\begin{remarkk}
{\rm{
\begin{enumerate}
\item
This result holds in fact for the Lusin spaces since they are
    Borel isomorphic to the Cantor set. Besides it extends easily to
    countable spaces. 
\item The particular case where $T=\calM_1(X)$, i.e., the space of
  probability measures on $X$ under the weak topology and $t\to \ga_t$
  being the identity map, is particularly
  important for the sequel. In this case we obtain a kernel $N$ such
  that  $(x,\ga)\to N_\ga f(x)$ is measurable and $N_\ga f$ is an
  $\calF$-measurable version of $E_\ga[f|\calF]$. 
\end{enumerate}
}}
\end{remarkk}

\begin{lemma}
\label{section-lemma}
Let $\rho$ and $\nu$ be two probability measures on $W$ such that
$$
d_H(\rho,\nu)<\infty
$$
and let $\ga\in\Sigma(\rho,\nu)$ be an optimal measure, i.e.,
$J(\ga)=d_H^2(\rho,\nu)$, where $J$ is given by 
(\ref{J-defn}). Assume that $F$ is a regular 
finite dimensional subspace of $W$ with the 
corresponding projection  $\pi_F$  from $W$ to $F$ and let
$\pi_F^\bot=I_W-\pi_F$ . Define  $p_F$ as  the 
projection from $W\times W$ onto $F$ with $p_F(x,y)=\pi_Fx$ and let
$p_F^\bot(x,y)=\pi_F^\bot x$.  
Consider the Borel disintegration 
\beaa
\ga(\cdot)&=&\int_{F^\bot\times W}\ga(\,\cdot|x^\bot)\ga^\bot(dz^\bot)\\
&=&\int_{F^\bot}\ga(\,\cdot|x^\bot)\rho^\bot(dx^\bot)
\eeaa
along the projection  of $W\times W$ on $F^\bot$, where
$\rho^\bot$ is the measure
$\pi_F^\bot\rho$,  $\ga(\cdot\,|x^\bot)$ denotes the regular 
conditional probability $\ga(\cdot\,|p_F^\bot =x^\bot)$ and $\ga^\bot$
is the measure $p_F^\bot\ga$. Then,
$\rho^\bot$ and $\ga^\bot$-almost surely  
$\ga(\,\cdot|x^\bot)$ is optimal on $(x^\bot+F)\times W$.
\end{lemma}
\proof
Let $p_1,\,p_2$ be the projections of $W\times W$ defined as
$p_1(x,y)=\pi_F(x)$ and $p_2(x,y)=\pi_F(y)$. 
Note first the  following obvious identity:
$$
p_1\ga(\cdot\,|x^\bot)=\rho(\cdot\,|x^\bot)\,,
$$
$\rho^\bot$ and $\ga^\bot$-almost surely.
Define the sets $B\subset F^\bot\times \calM_1(F\times F)$ and $C$ as 
\beaa
B&=&\{(x^\bot,\theta):\,\theta\in\Sigma(p_1\ga(\cdot\,|x^\bot),
p_2\ga(\cdot\,|x^\bot))\}\\
C&=&\{(x^\bot,\theta)\in B:\,J(\theta)<J(\ga(\cdot\,|x^\bot)\}\,,
\eeaa
where $\calM_1(F\times F)$ denotes the set of probability measures on
$F\times F$. Let $K$ be the projection of $C$ on $F^\bot$. Since $B$
and $C$ are Borel measurable, $K$is a
Souslin set, hence it is $\rho^\bot$-measurable.  The selection
theorem (cf. \cite{Del-M}) implies the existence of a measurable map
$$
x^\bot\to \theta_{x^\bot}
$$ 
from $K$ to $\calM_1(F\times F)$ such that, $\rho^\bot$-almost
surely,  $(x^\bot,\theta_{x^\bot})\in C$. Define 
$$
\theta(\cdot)=\int_K\theta_{x^\bot}(\cdot)d\rho^\bot(x^\bot)+
\int_{K^c}\ga(\cdot\,|x^\bot)d\rho^\bot(x^\bot)\,.
$$
Then $\theta\in \Sigma(\rho,\nu)$ and we have 
\beaa
J(\theta)&=&\int_K J(\theta_{x^\bot})d\rho^\bot(x^\bot)
+\int_{K^c} J(\ga(\cdot\,|x^\bot))d\rho^\bot(x^\bot)\\
&<&\int_{K} J(\ga(\cdot\,|x^\bot))d\rho^\bot(x^\bot)
+\int_{K^c} J(\ga(\cdot\,|x^\bot))d\rho^\bot(x^\bot)\\
&=&J(\ga)\,,
\eeaa
hence we obtain $J(\theta)<J(\ga)$ which is a contradiction to the
optimality of $\ga$.
\qed

\begin{lemma}
\label{f-dim-lemma}
Assume that  the hypothesis of Lemma \ref{section-lemma} holds  and 
let $F$ be  any regular  finite dimensional  subspace  of $W$. 
Denote by $\pi_F$ the projection operator associated to it and let
$\pi_F^\bot=I_W-\pi_F$. If $\pi_F^\bot\rho$-almost surely,  the
regular conditional probability 
$\rho(\cdot\,|\pi_F^\bot=x^\bot )$ vanishes on the subsets of
$x^\bot+F$ whose  Hausdorff dimension are at most equal to 
${\mbox{\rm  dim}}(F)-1$,  then
there exists a map $T_F:F\times F^\bot\to F$ such that 
$$
\ga\left(\left\{(x,y)\in W\times
  W:\,\pi_F y=T_F(\pi_Fx,\pi_F^\bot x)\right\}\right)=1\,.
$$
\end{lemma}
\proof
Let $C_{x^\bot}$ be the support of the regular conditional probability
$\ga(\cdot\,|x^\bot)$ in $(x^\bot+F)\times W$. 
We know from Lemma \ref{section-lemma}
that the measure  $\ga(\cdot\,|x^\bot)$ is optimal in
$\Sigma(\pi_1\ga(\cdot\,|x^\bot),\pi_2\ga(\cdot\,|x^\bot))$, 
with $J(\ga(\cdot\,|x^\bot))<\infty$ for $\rho^\bot$-almost everywhere
$x^\bot$. From  Theorem 2.3 of \cite{G-Mc} and from \cite{Ab-H}, the
set $C_{x^\bot}$ 
is cyclically monotone, moreover, $C_{x^\bot}$ is a subset of
$(x^\bot+F)\times H$, 
hence the cyclic monotonicity of it implies that the  set
$K_{x^\bot}\subset F\times F$,  defined as 
$$
K_{x^\bot}=\{(u,\pi_Fv)\in F\times F:\,(x^\bot+u,v)\in C_{x^\bot}\}
$$
is cyclically monotone in $F\times F$. Therefore $K_{x^\bot}$ is 
included in the subdifferential of a convex function defined on $F$. 
Since, by hypothesis,  the first marginal of $\ga(\cdot\,|x^\bot)$, i.e.,
$\rho(\cdot\,|x^\bot)$ vanishes on the subsets of $x^\bot+F$ 
of co-dimension one,  the subdifferential
under question, denoted as $U_F(u,x^\bot)$  is
$\rho(\cdot\,|x^\bot)$-almost surely univalent (cf. \cite{An-K,Mc1}). This
implies that  
$$
\ga(\cdot\,|x^\bot)\left(\left\{(u,v)\in
    C_{x^\bot}:\,\pi_Fv=U_F(u,x^\bot)\right\}\right)=1\,, 
$$
$\rho^\bot$-almost surely.
Let 
$$
K_{x^\bot,u}=\left\{v\in W:\,(u,v)\in K_{x^\bot}\right\}\,.
$$
Then $K_{x^\bot,u}$ consists of a single point for
almost all $u$ with respect to $\rho(\cdot\,|x^\bot)$.
Let 
$$
N=\left\{(u,x^\bot)\in F\times
  F^\bot:\,{\mbox{Card}}(K_{x^\bot,u})>1\right\}\,,
$$
note that $N$ is a Souslin set, hence it is universally
measurable. Let $\sigma$ be the measure which is defined as the image
of $\rho$ under the projection $x\to (\pi_Fx,\pi_F^\bot x)$. We then
have 
\beaa
\sigma(N)&=&\int_{F^\bot}\rho^\bot(dx^\bot)\int_F
\won_N(u,x^\bot)\rho(du|x^\bot)\\
&=&0\,.
\eeaa
Hence $(u,x^\bot)\mapsto K_{x^\bot,u}=\{y\}$ is $\rho$ and  $\ga$-almost surely
well-defined and it suffices to denote this map by $T_F$ to achive the
proof.
\qed

\begin{theorem}
\label{monge-general}
Suppose that $\rho$ and $\nu$ are two probability measures on
$W$ such  that
$$
d_H(\rho,\nu)<\infty\,.
$$
Let $(\pi_n,n\geq 1)$ be a total increasing  sequence of regular
projections (of $H$, converging to the identity map of $H$). 
Suppose  that, for any $n\geq 1$, the regular
conditional probabilities $\rho(\cdot\,|\pi_n^\bot=x^\bot)$ vanish
$\pi_n^\bot\rho$-almost surely on 
the subsets of  $(\pi_n^\bot)^{-1}(W)$ with Hausdorff dimension
$n-1$. Then there exists a    
unique solution of the   Monge-Kantorovitch problem, denoted by $\ga\in
\Sigma(\rho,\nu)$ and  $\ga$ is supported by the graph of a Borel
map $T$ which is the solution of the
Monge problem.  $T:W\to W$ is  of the form $T=I_W+\xi$ , where $\xi\in
H$ almost surely. Besides  we have 
\beaa
d_H^2(\rho,\nu)&=&\int_{W\times W}|T(x)-x|_H^2d\ga(x,y)\\
&=&\int_{W}|T(x)-x|_H^2d\rho(x)\,, 
\eeaa
and  
for $\pi_n^\bot\rho$-almost almost all $x_n^\bot$, the map $u\to
\xi(u+x_n^\bot)$ is cyclically monotone on
$(\pi_n^\bot)^{-1}\{x_n^\bot\}$, in the sense that 
$$
\sum_{i=1}^N\left(\xi(x_n^\bot+u_i),u_{i+1}-u_i\right)_H\leq 0
$$
$\pi_n^\bot\rho$-almost surely, for any cyclic sequence
$\{u_1,\ldots,u_N,u_{N+1}=u_1\}$ from $\pi_n(W)$. Finally, if, for any $n\geq
1$, $\pi_n^\bot\nu$-almost surely,  $\nu(\cdot\,|\pi_n^\bot=y^\bot)$
 also vanishes on the $n-1$-Hausdorff dimensional  subsets  
 of $(\pi_n^\bot)^{-1}(W)$, then $T$ is invertible, i.e, there exists
 $S:W\to W$ of the form $S=I_W+\eta$ such that  $\eta\in H$ satisfies
 a similar  cyclic monotononicity property as $\xi$ and that 
\beaa
1&=&\ga\left\{(x,y)\in W\times W: T\circ S(y)=y\right\}\\
&=&\ga\left\{(x,y)\in W\times W: S\circ T(x)=x\right\}\,.
\eeaa
In particular we have 
\beaa
d_H^2(\rho,\nu)&=&\int_{W\times W}|S(y)-y|_H^2d\ga(x,y)\\
&=&\int_{W}|S(y)-y|_H^2d\nu(y)\,. 
\eeaa
\end{theorem}
\begin{remarkk}
In particular, for all  the measures $\rho$ which are  absolutely
continuous with respect to the  Wiener measure $\mu$,  the second
hypothesis is satisfied, i.e., the measure
$\rho(\cdot\,|\pi_n^\bot=x_n^\bot)$ vanishes on the sets of Hausdorff
dimension $n-1$.
\end{remarkk} 
\proof
Let $(F_n,n\geq 1)$ be the  increasing sequence of regular subspaces 
associated to $(\pi_n,n\geq 1)$,  
 whose union is dense in $W$.  From Lemma \ref{f-dim-lemma}, for any
 $F_n$, there exists a map $T_n$, such that
 $\pi_n y=T_n(\pi_n x,\pi_n^\bot x)$ for $\ga$-almost all $(x,y)$, where
 $\pi_n^\bot=I_W-\pi_n$. Write 
 $T_n$ as $I_n+\xi_n$, where $I_n$ denotes the identity map on
 $F_n$. Then we have the following representation:
$$
\pi_ny=\pi_nx+\xi_n(\pi_nx,\pi_n^\bot x)\,,
$$
$\ga$-almost surely.
Since 
\beaa
\pi_n y-\pi_n x&=&\pi_n( y-x)\\
&=&\xi_n(\pi_n x,\pi_n^\bot x)\,
\eeaa 
and since $y-x\in H$ $\ga$-almost surely, $(\pi_ny-\pi_nx,n\geq 1)$
converges $\ga$-almost surely.
Consequently $(\xi_n,n\geq 1)$ converges $\ga$, hence $\rho$ almost
surely to a measurable  $\xi$. Consequently we obtain 
$$
\ga\left(\left\{(x,y)\in W\times W:\, y=x+\xi(x)\right\}\right)=1\,.
$$
Since $J(\ga)<\infty$, $\xi$ takes its values almost surely in the
Cameron-Martin space $H$. The cyclic monotonicity of $\xi$  is
obvious. To prove the uniqueness, assume that we have two 
optimal solutions $\ga_1$ and $\ga_2$ with the same marginals and
$J(\ga_1)=J(\ga_2)$. Since $\beta\to J(\beta)$ is linear,  the measure
defined as  $\ga=\frac{1}{2}(\ga_1+\ga_2)$ is
also optimal and it has also the same marginals $\rho$ and
$\nu$. Consequently, it is also  supported by the graph of a map 
$T$. Note that  $\ga_1$ and $\ga_2$ are absolutely continuous with
respect to $\ga$,  let $L_1(x,y)$ be the Radon-Nikodym density of
$\ga_1$ with respect to $\ga$. For any $f\in C_b(W)$, we then  have 
\beaa
\int_Wfd\rho&=&\int_{W\times W}f(x)d\ga_1(x,y)\\
&=&\int_{W\times W}f(x)L_1(x,y)d\ga(x,y)\\
&=&\int_Wf(x)L_1(x,T(x))d\rho(x)\,. 
\eeaa
Therefore we should have $\rho$-almost surely,  $L_1(x,T(x))=1$, hence
also $L_1=1$ almost everywhere $\ga$ and this implies that
$\ga=\ga_1=\ga_2$. The second part about the invertibility of $T$ is
totally symmetric, hence its proof follows along the same lines as the
proof for $T$.
\qed

\begin{corollary}
\label{c-mono}
Assume that $\rho$ is equivalent to the Wiener measure $\mu$, then for
any $h_1,\ldots,h_N\in H$ and for any permutation $\tau$ of
$\{1,\ldots,N\}$, we have, with the notations of Theorem
\ref{monge-general}, 
$$
\sum_{i=1}^N\left(h_i+\xi(x+h_i),h_{\tau(i)}-h_i\right)_H\leq 0
$$
$\rho$-almost surely.
\end{corollary}
\proof
Again with the notations of the theorem, $\rho_k^\bot$-almost surely,
the graph of the map $x_k\to x_k+\xi_k(x_k,x_k^\bot)$ is cyclically
monotone on $F_k$. Hence, for the case $h_i\in F_n$ for all
$i=1,\ldots,N$ and $n\leq k$, we have 
$$
\sum_{i=1}^N\left(h_i+x_k+\xi_k(x_k+h_i,x_k^\bot),h_{\tau(i)}-h_i\right)_H\leq 0\,.
$$
Since $\sum_i(x_k,h_{\tau(i)}-h_i)_H=0$, we also have 
$$
\sum_{i=1}^N\left(h_i+\xi_k(x_k+h_i,x_k^\bot),h_{\tau(i)}-h_i\right)_H\leq 0\,.
$$
We know that $\xi_k(x_k+h_i,x_k^\bot)$ converges to $\xi(x+h_i)$
$\rho$-almost surely. Moreover $h\to \xi(x+h)$ is continuous from $H$
to $L^0(\rho)$ and the proof follows.
\qed


\section{The Monge-Amp\`ere equation}
\label{equation}
\label{monge-ampere}


Assume   that $W=\R^n$ and take a   density $L\in
\LL\log\LL$. Let $\phi\in \DD_{2,1}$ be the $1$-convex function such
that $T=I+\nabla \phi$ maps $\mu$ to $L\cdot \mu$. Let $S=I+\nabla\psi$
be its inverse with $\psi\in \DD_{2,1}$. 
Let now $\nabla_a^2\phi$ be the
second Alexandrov derivative of 
$\phi$, i.e.,  the Radon-Nikodym derivative of the absolutely
continuous part of the 
vector measure $\nabla^2\phi$ with respect to the Gaussian measure
$\mu$ on $\reals^n$. Since $\phi$ is $1$-convex, it follows that
$\nabla^2\phi\geq -I_{\reals^n}$ in the sense of the distributions,
consequently $\nabla_a^2\phi\geq -I_{\reals^n}$ $\mu$-almost
surely. Define also the Alexandrov version $\calL_a\phi$ of
$\calL\phi$ as the Radon-Nikodym derivative of  the absolutely
continuous part of the distribution 
$\calL\phi$. Since we are in finite dimensional situation, we have the
explicit expression for $\calL_a\phi$ as 
$$
\calL_a\phi(x)=(\nabla\phi(x),x)_{\reals^n}-{\rm
  trace}\left(\nabla_a^2\phi\right)\,.
$$
Let  $\La$ be the Gaussian  Jacobian 
$$
\La=\dett\left(I_{\R^n}+\nabla^2_a\phi\right)\exp\left\{-\calL_a\phi-
\frac{1}{2}|\nabla\phi|_{\R^n}^2\right\}\,.
$$
\begin{remarkk}{\rm
In this expression as well as in the sequel,  the notation
$\dett(I_H+A)$ denotes the modified 
Carleman-Fredholm determinant of the operator $I_H+A$ on a Hilbert
space $H$. If $A$ is an operator of finite rank, then it  is defined  as 
$$
\dett\left(I_H+A\right)=\prod_{i=1}^n(1+l_i)e^{-l_i}\,,
$$
where $(l_i,\,i\leq n)$ denotes the eigenvalues of $A$ counted with
respect to their multiplicity. In fact this determinant has an
analytic extension to the space of Hilbert-Schmidt operators on a
separable Hilbert space, cf.  \cite{Du-S} and Appendix A.2 of
\cite{BOOK}. As explained in \cite{BOOK}, the modified determinant
exists for the Hilbert-Schmidt operators while the ordinary
determinant does not, since the latter  requires the existence of the trace of
$A$. Hence the modified Carleman-Fredholm determinant  is particularly
useful when one studies the absolute 
continuity properties of the image of a Gaussian measure under
non-linear transformations in the setting of infinite dimensional
Banach spaces (cf., \cite{BOOK} for further information).  }
\end{remarkk}
It follows from the change of variables formula given in  Corollary
4.3  of \cite{Mc2}, that, for any $f\in C_b(\R^n)$,  
$$
E[f\circ T \,\La]=E\left[f\,1_{\partial\Phi(M)}\right]\,,
$$
where $M$ is the set of non-degeneracy of $I_{\R^n}+\nabla_a^2\phi$,
$$
\Phi(x)=\frac{1}{2}|x|^2+\phi(x)
$$ 
and $\partial \Phi$ denotes  the subdifferential of the convex
function  $\Phi$. Let us note that, in case $L>0$ almost surely, $T$
has  a global inverse $S$, i.e., $S\circ T=T\circ S=I_{\R^n}$
$\mu$-almost surely and $\mu(\partial\Phi(M))=\mu(S^{-1}(M))$. 
Assume now that $\La>0$ almost surely, i.e., that $\mu(M)=1$. Then,
for any $f\in C_b(\R^n)$, we have 
\beaa
E[f\circ T]&=&E\left[f\circ T\,\frac{\La}{\La\circ T^{-1}\circ T}\right]\\
&=&E\left[f\,\frac{1}{\La\circ T^{-1}}1_{\partial\Phi(M)}\right]\\
&=&E[f\,L]\,,
\eeaa
where $T^{-1}$ denotes the left inverse of $T$ whose existence is
guaranteed by Theorem \ref{gaussian-case}. Since $T(x)\in
\partial\Phi(M)$ almost surely, it follows from the above calculations
$$
\frac{1}{\La}=L\circ T\,,
$$
almost surely. Take now  any $t\in [0,1)$, the map
$x\to \frac{1}{2}|x|_H^2+t\phi(x)=\Phi_t(x)$ is strictly convex and a simple
calculation implies that the mapping  $T_t=I+t\nabla\phi$ is
$(1-t)$-monotone (cf. \cite{BOOK}, Chapter 6), consequently it has a
left inverse denoted by $S_t$. Let us  denote by $\Psi_t$ the  Legendre
transformation of $\Phi_t$: 
$$
\Psi_t(y)=\sup_{x\in \R^n}\left\{(x,y)-\Phi_t(x)\right\}\,.
$$
A simple calculation shows that 
\beaa
\Psi_t(y)&=&\sup_x\left[(1-t)\left\{(x,y)-\frac{|x|^2}{2}\right\}+
t\left\{(x,y)-\frac{|x|^2}{2}-\phi(x)\right\}\right]\\
&\leq&(1-t)\frac{|y|^2}{2}+t\Psi_1(y)\,.
\eeaa
Since $\Psi_1$ is the Legendre transformation of
$\Phi_1(x)=|x|^2/2+\phi(x)$ and since $L\in \LL\log \LL$, it is
finite on a convex set of full measure, hence it is finite everywhere.
Consequently $\Psi_t(y)<\infty$ for any $y\in \R^n$. Since a finite,
convex function is almost everywhere differentiable, $\nabla \Psi_t$
exists almost everywhere on  and it is equal almost everywhere on
$T_t(M_t)$ to the left inverse $T_t^{-1}$, where $M_t$ is the set of
non-degeneracy of $I_{\R^n}+t\nabla^2_a\phi$. Note that $\mu(M_t)=1$. 
 The strict convexity implies that $T_t^{-1}$ is Lipschitz with a
Lipschitz constant $\frac{1}{1-t}$. 
Let now $\La_t$ be the Gaussian  Jacobian 
$$
\La_t=\dett\left(I_{\R^n}+t\nabla^2_a\phi\right)\exp\left\{-t\calL_a\phi-
\frac{t^2}{2}|\nabla\phi|_{\R^n}^2\right\}\,.
$$
Since the domain of $\phi$ is the whole space $\R^n$, $\La_t>0$ almost
surely, hence, as we have explained above,   it follows from the
change of variables formula of \cite{Mc2} that $T_t\mu$ is absolutely
continuous with respect to $\mu$ and that 
$$
\frac{1}{\La_t}=L_t\circ T_t\,,
$$
$\mu$-almost surely. 

\noindent
Let us come back to the infinite dimensional case: we first  give an
inequality which may be useful.
\begin{theorem}
\label{monge1-thm}
Assume that $(W,\mu,H)$ is an abstract  Wiener space, assume
that $K,L\in \LL_+^1(\mu)$ with $K>0$ almost surely  and denote by
$T:W\to W$ the transfer map 
$T=I_W+\nabla \phi$, which maps the measure $Kd\mu$ to the measure
$Ld\mu$. Then the following inequality holds:
\begin{equation}
\label{monge1-ineq}
\frac{1}{2}E[|\nabla\phi|_H^2]\leq E[-\log K+\log L\circ T]\,.
\end{equation} 
\end{theorem}
\proof
Let us define $k$ as $k=K\circ T^{-1}$, then for any $f\in C_b(W)$, we
have 
\beaa
\int_Wf(y)L(y)d\mu(y)&=&\int_Wf\circ T(x)K(x)d\mu(x)\\
&=&\int_Wf\circ T(x)k\circ T(x)d\mu(x)\,,
\eeaa
hence 
$$
T\mu=\frac{L}{k}\,.\mu\,.
$$
It then follows from the inequality \ref{tal-ineq} that 
\beaa
\frac{1}{2}E\left[|\nabla\phi|_H^2\right]&\leq&E\left[\frac{L}{k}\log\frac{L}{k}\right]\\
&=&E\left[\log\frac{L\circ T}{k\circ T}\right]\\
&=&E[-\log K+\log L\circ T]\,.
\eeaa
\qed

\noindent
 Suppose that $\phi\in \DD_{2,1}$ is a $1$-convex Wiener
 functional. Let $V_n$ be 
the sigma algebra generated by $\{\delta e_1,\ldots,\delta e_n\}$,
where $(e_n,\,n\geq 1)$ is an orthonormal basis of the Cameron-Martin
space $H$. Then $\phi_n=E[\phi|V_n]$ is again $1$-convex
(cf.\cite{F-U1}), hence $\calL\phi_n$ is a measure as it can be easily
verified.  However the sequence
$(\calL\phi_n,\,n\geq 1)$ converges to $\calL\phi$ only  in
$\DD'$. Consequently, there is no reason for the limit 
$\calL\phi$ to be a  measure. In case this happens, we shall denote the
Radon-Nikodym density with respect to $\mu$,  of the 
absolutely continuous part  of this measure  by $\calL_a\phi$. 
\begin{lemma}
\label{app-lemma}
Let $\phi\in \DD_{2,1}$ be $1$-convex and let $V_n$ be defined as
above  and define $F_n=E[\phi|V_n]$. Then
the sequence $(\calL_a F_n,n\geq 1)$ 
is a submartingale, where $\calL_aF_n$ denotes the
$\mu$-absolutely continuous part of the measure $\calL F_n$.
\end{lemma}
\proof
Note that, due to the $1$-convexity, we have $\calL_aF_n\geq
\calL F_n$ for any $n\in \NN$. 
Let $X_n=\calL_aF_n$ and $f\in \DD$ be a positive, $V_n$-measurable
test function. 
Since $\calL E[\phi|V_n]=E[\calL\phi|V_n]$, we have 
\beaa
E[X_{n+1}\,f]&\geq&\langle \calL F_{n+1},f\rangle\\
&=&\langle \calL F_{n},f\rangle\,,
\eeaa
where $\langle\cdot,\cdot\rangle$ denotes the duality bracket for the
dual pair $(\DD',\DD)$. 
Consequently  
$$
E[f\,E[X_{n+1}|V_n]]\geq \langle \calL F_{n},f\rangle\,,
$$ 
for any positive, $V_n$-measurable test function $f$, it follows that  the
absolutely continuous part of $\calL F_n$  is also  dominated by the
same conditional expectation and  this proves the submartingale
property. 
\qed

\begin{lemma}
\label{Lf-lemma}
Assume that $L\in \LL\log \LL$ is a positive random variable whose
expectation is one. Assume further that it  is lower bounded by a  constant
$a>0$. Let $T=I_W+\nabla \phi$ be the transport map such that
$T\mu=L\,.\mu$ and let $T^{-1}=I_W+\nabla \psi$. Then $\calL \psi$ is
a Radon measure on $(W,\calB(W))$. If $L$ is upper bounded by $b>0$,
then $\calL\phi$ is also a Radon measure on $(W,\calB(W))$.
\end{lemma}
\proof
Let $L_n=E[L|V_n]$, then $L_n\geq a$ almost surely. Let
$T_n=I_W+\nabla\phi_n$ be the transport map which satisfies
$T_n\mu=L_n\,.\mu$ and let $T_n^{-1}=I_W+\nabla\psi_n$ be its
inverse. We have 
$$
L_n=\dett\left(I_H+\nabla_a^2\psi_n\right)\exp\left[-\calL_a\psi_n-\frac{1}{2}|\nabla\psi_n|_H^2\right]\,.
$$
By the hypothesis $-\log L_n\leq -\log a$.  Since $\psi_n$ is
$1$-convex, it follows from the finite dimensional results that
$\dett\left(I_H+\nabla_a^2\psi_n\right)\in [0,1]$ almost surely. 
Therefore  we have
$$
\calL_a\psi_n\leq -\log a\,,
$$
besides $\calL\psi_n\leq \calL_a\psi_n$ as distributions, consequently
$$
\calL\psi_n\leq -\log a
$$
as distributions, for any $n\geq 1$. Since $\lim_n\calL\psi_n=\calL
\psi$ in $\DD'$, we obtain $\calL\psi\leq -\log a$, hence $-\log
a-\calL\psi\geq 0$ as a distribution, hence $\calL\psi$  is a Radon measure on
$W$, c.f., \cite{F-P}, \cite{ASU}.  This proves the
first claim. Note that whenever $L$ is upperbounded, $\La=1/L\circ T$
is lowerbounded, hence  the proof of the second claim  is similar to
that of the first one.
\qed

\begin{theorem}
\label{girsanov-thm}
Assume that $L$ is a strictly  positive bounded  random variable with
$E[L]=1$. Let $\phi\in \DD_{2,1}$ be
the $1$-convex Wiener functional such that 
$$
T=I_W+\nabla \phi
$$ 
is the transport map realizing the measure $L\,.\mu$ and let
$S=I_W+\nabla\psi$ be 
its inverse. Define $F_n=E[\phi|V_n]$, then the submartingale
$(\calL_aF_n,n\geq 1)$ converges almost surely to $\calL_a\phi$. Let
$\la(\phi)$ be  the random variable defined as 
\beaa
\la(\phi)&=&\lim\inf_{n\to\infty}\La_n\\
&=&\left(\lim\inf_n\dett\left(I_H+\nabla^2_aF_n\right)\right)\exp\left\{-\calL_a\phi-\half|\nabla\phi|_H^2\right\}
\eeaa
where
$$
\La_n=\dett\left(I_H+\nabla_a^2F_n\right)\exp\left\{-\calL_aF_n-\half|\nabla
  F_n|_H^2\right\}\,.
$$
Then it holds true that 
\begin{equation}
\label{sub-gir}
E[f\circ T\,\la(\phi)]\leq E[f]
\end{equation}
for any $f\in C^+_b(W)$, in particular $\la(\phi)\leq \frac{1}{L\circ
  T}$ almost surely. If $E[\la(\phi)]=1$, then the inequality in
(\ref{sub-gir}) becomes an equality and  we also  have  
$$
\la(\phi)=\frac{1}{L\circ T}\,.
$$
\end{theorem}
\proof Let us remark that, due to the $1$-convexity,
$0\leq\dett\left(I_H+\nabla^2_aF_n\right)\leq 1$, hence the $\lim\inf$
exists. Now,  Lemma \ref{Lf-lemma} implies that $\calL\phi$ is a Radon
measure. Let $F_n=E[\phi|V_n]$, then we know from Lemma
\ref{app-lemma} that $(\calL_a F_n,n\geq 1)$ is a submartingale. Let
$\calL^+\phi$ denote the positive part of the measure $\calL\phi$. Since
$\calL^+\phi\geq \calL\phi$, we have also $E[\calL^+\phi|V_n]\geq
E[\calL \phi|V_n]=\calL F_n$. This implies that
$E[\calL^+\phi|V_n]\geq\calL_a^+F_n$. Hence we find that 
$$
\sup_n E[\calL^+_a F_n]<\infty
$$
and this condition implies that the submartingale $(\calL_aF_n,n\geq 1)$
converges almost surely. We shall now identify the limit of this
submartingale. Let $\calL_s G$ be  the singular part of the
measure $\calL G$ for a  Wiener function $G$ such that $\calL G$ is a
measure. We have
\beaa
E[\calL\phi|V_n]&=&E[\calL_a\phi|V_n]+E[\calL_s\phi|V_n]\\
&=&\calL_aF_n+\calL_sF_n\,,
\eeaa
hence 
$$
\calL_aF_n=E[\calL_a\phi|V_n]+E[\calL_s\phi|V_n]_a
$$
almost surely, where $E[\calL_s\phi|V_n]_a$ denotes the absolutely
continuous part of the measure $E[\calL_s\phi|V_n]$. Note that, from
the Theorem of Jessen 
(cf., for example Theorem 1.2.1 of \cite{BOOK}),
$\lim_nE[\calL_s^+\phi|V_n]_a=0$ and
$\lim_nE[\calL_s^-\phi|V_n]_a=0$ almost surely, hence we have 
$$
\lim_n\calL_aF_n=\calL_a\phi\,,
$$
$\mu$-almost surely. To complete the proof, an application of the
Fatou lemma implies that  
\beaa
E[f\circ T\,\la(\phi)]&\leq&E[f]\\
&=&E\left[f\circ T\,\frac{1}{L\circ T}\right]\,,
\eeaa
for any $f\in C_b^+(W)$. Since $T$ is invertible, it follows that  
$$
\la(\phi)\leq \frac{1}{L\circ T}
$$
almost surely. Therefore, in  case $E[\la(\phi)]=1$,  we have
$$
\la(\phi)=\frac{1}{L\circ T}\,,
$$
and this completes the proof.
\qed

\begin{corollary}
\label{monge-amp.eqn}
Assume that $K,L$ are two positive random variables with values in a
bounded interval $[a,b]\subset (0,\infty)$ such that
$E[K]=E[L]=1$. Let
$T=I_W+\nabla\phi$, $\phi\in \DD_{2,1}$,  be the
transport map pushing $Kd\mu$ to $Ld\mu$, i.e, $T(Kd\mu)=Ld\mu$. We
then have
$$
L\circ T\,\la(\phi)\leq K\,,
$$
$\mu$-almost surely. In particular, if $E[\la(\phi)]=1$, then $T$ is
the solution of the Monge-Amp\`ere equation.
\end{corollary}
\proof Since $a>0$, 
$$
\frac{dT\mu}{d\mu}=\frac{L}{K\circ T}\leq \frac{b}{a}\,.
$$
Hence,  Theorem \ref{sub-gir} implies that 
\beaa
E[f\circ T\,L\circ T\,\la(\phi)]&\leq& E[f\,L]\\
&=&E[f\circ T\,K]\,,
\eeaa
consequently
$$
L\circ T\,\la(\phi)\leq K\,,
$$
the rest of the claim is now obvious.
\qed

For later use we give also the folowing result:
\begin{theorem}
\label{interpol-thm}
Assume that $L$ is a positive random variable of class $\LL\log \LL$
such that $E[L]=1$. Let $\phi\in \DD_{2,1}$ be the $1$-convex function
corresponding to the transport map $T=I_W+\nabla\phi$. Define
$T_t=I_W+t\nabla\phi$, where  $t\in [0,1]$. Then,
for any $t\in [0,1]$, $T_t\mu$ is absolutely continuous with respect
to the Wiener measure $\mu$.
\end{theorem}
\proof
Let $\phi_n$ be defined as the transport map corresponding to
$L_n=E[P_{1/n}L_n|V_n]$ and define $T_n$ as $I_W+\nabla\phi_n$. For $t\in
[0,1)$, let $T_{n,t}=I_W+t\nabla\phi_n$.  It follows from the finite
dimensional results which are summarized  in the
beginning of this section, that $T_{n,t}\mu$ is absolutely continuous 
with respect to $\mu$. Let $L_{n,t}$ be the corresponding
Radon-Nikodym density and define $\La_{n,t}$ as 
$$
\La_{n,t}=\dett\left(I_H+t\nabla^2_a\phi_{n}\right)
\exp\left\{-t\calL_a\phi_n-\frac{t^2}{2}|\nabla\phi_n|_H^2\right\}\,.
$$
Besides, for any $t\in [0,1)$, 
\begin{equation}
\label{mon-con}
\left((I_H+t\nabla_a^2\phi_n)h,h\right)_H>0\,,
\end{equation}
$\mu$-almost surely for any $0\neq h\in H$. Since $\phi_n$ is of finite 
rank, \ref{mon-con} implies that $\La_{n,t}>0$ $\mu$-almost surely
and we have shown  at the beginning of this section 
$$
\La_{n,t}=\frac{1}{L_{n,t}\circ T_{n,t}}
$$
$\mu$-almost surely.  
An easy calculation shows that $t\to \log\dett(I+t\nabla_a^2\phi_n)$ is
a non-increasing function. Since $\calL_a\phi_n\geq \calL\phi_n$, we
have $E[\calL_a\phi_n]\geq 0$. Consequently
\beaa
E\left[L_{t,n}\log L_{t,n}\right]&=&E\left[\log L_{n,t}\circ T_{n,t}\right]\\
&=&-E\left[\log\La_{t,n}\right]\\
&=&E\left[-\log\dett\left(I_H+t\nabla^2\phi_n\right)+t\calL_a\phi_n
+\frac{t^2}{2}|\nabla\phi_n|_H^2\right]\\
&\leq&E\left[-\log\dett\left(I_H+\nabla^2\phi_n\right)+\calL_a\phi_n
+\frac{1}{2}|\nabla\phi_n|_H^2\right]\\
&=&E\left[L_n\log L_n\right]\\
&\leq&E[L\log L]\,,
\eeaa
by the  Jensen inequality. Therefore 
$$
\sup_nE[L_{n,t}\log L_{n,t}]<\infty
$$
and this implies that the sequence $(L_{n,t},n\geq 1)$ is uniformly
integrable for any $t\in [0,1]$. Consequently it has a subsequence
which converges weakly in $L^1(\mu)$ to some $L_t$.  Since, from Theorem
\ref{gaussian-case},  $\lim_n\phi_n=\phi$ in $\DD_{2,1}$, where $\phi$
is the transport map associated to $L$, for any $f\in C_b(W)$, we have
\beaa
E[f\circ T_t]&=&\lim_k E\left[f\circ T_{n_k,t}\right]\\
&=&\lim_k E\left[f\,L_{n_k,t}\right]\\
&=&E[f\,L_t]\,,
\eeaa
hence the theorem is proved.
\qed

\subsection{The solution of the Monge-Amp\`ere equation via
  Ito-renormalization} 
\label{sub-equation}
We can interpret the Monge-Amp\`ere equation as follows: given two
probability densities $K$ and $L$, find a map $T:W\to W$ such that 
$$
L\circ T\, J(T)=K
$$
almost surely, where $J(T)$ is a kind of Jacobian to be written in
terms of $T$. In Corollary \ref{monge-amp.eqn}, we have shown the
existence of some $\la(\phi)$ which gives an inequality instead of the
equality. Although in the finite dimensional case there are some
regularity results about the transport map (cf., \cite{Caf}), in the
infinite dimensional case  such techniques do not work.  All these
difficulties  can
be circumvented using the miraculous 
renormalization of the Ito calculus. In fact assume that $K$ and $L$
satisfy the hypothesis of the corollary. First let us indicate that we
can assume $W=C_0([0,1],\reals)$ (cf., \cite{BOOK}, Chapter II, to see
how one can pass from an abstract Wiener space to the standard
one) and in this case the Cameron-Martin space $H$ becomes
$H^1([0,1])$, which is the space of absolutely continuous functions on
$[0,1]$, with a square integrable Sobolev derivative. Let now 
$$
\La=\frac{K}{L\circ T}\,,
$$
where $T$ is as  constructed above. Then $\La.\mu$ is a Girsanov
measure for the map $T$. This means that the law of the stochastic
process $(t,x)\to T_t(x)$ under $\La.\mu$ is equal to the Wiener
measure, where $T_t(x)$ is defined as the evaluation of the trajectory
$T(x)$ at $t\in [0,1]$. In other words the process $(t,x)\to T_t(x)$
is a Brownian motion under the probability $\La.\mu$. Let
$(\calF^T_t,t\in [0,1])$ be its filtration, the invertibility of  $T$
implies  that 
$$
\bigvee_{t\in [0,1]}\calF^T_t=\calB(W)\,.
$$
 $\La$ is upper and lower  bounded $\mu$-almost surely, hence also
 $\La.\mu$-almost surely. The Ito representation theorem implies that
 it can be represented as 
$$
\La=E[\La^2]\exp\left\{-\int_0^1\dot{\alpha}_sdT_s-\half\int_0^1|\dot{\alpha}_s|^2ds\right\}\,,
$$
where $\alpha(\cdot)=\int_0^\cdot\dot{\alpha}_sds$ is an $H$-valued
random variable. In fact $\alpha$ can be calculated explicitly using
the Ito-Clark representation theorem (cf., \cite{ASU}), and it is
given as 
\begin{equation}
\label{alpha}
\dot{\alpha}_t=\frac{E_\La[D_t\La|\calF^T_t]}{E_\La[\La|\calF^T_t]}
\end{equation}
$dt\times \La d\mu$-almost surely, where $E_\La$ denotes the
expectation operator with respect to $\La.\mu$ and $D_t\La$ is the
Lebesgue density of the absolutely continuous map $t\to
\nabla\La(t,x)$. From the relation (\ref{alpha}), it follows that  $\alpha$ is 
a function of $T$, hence we have obtained the strong solution of the
Monge-Amp\`ere equation. Let us announce all this as 
\begin{theorem}
\label{monge-amp.2}
Assume that $K$ and $L$ are upper and lower bounded densities, let $T$
be the transport map constructed in Theorem \ref{monge-general}. Then
$T$ is also the strong solution of the Monge-Amp\`ere equation in the
Ito sense, namely 
$$
E[\La^2]\,L\circ T
\exp\left\{-\int_0^1\dot{\alpha}_sdT_s-\half\int_0^1|\dot{\alpha}_s|^2ds\right\}=K\,,
$$
$\mu$-almost surely, where $\alpha$ is given with (\ref{alpha}).
\end{theorem}

\noindent
{\bf{Acknowledgement:}} The authors are grateful to Fran\c{c}oise
Combelles  for all the bibliographical help that she has supplied for
the realization of this research.

\noindent{\small{
\begin{itemize}
\item D. Feyel, Universit\'e d'Evry-Val-d'Essone, 91025 Evry Cedex, France.
E-mail: feyel@maths.univ-evry.fr
\item A. S. \"Ust\"unel, ENST, D\'ept. Infres, 46, rue Barrault, 75634
  Paris Cedex 13, France.
E-mail: ustunel@enst.fr
\end{itemize}
}}


\begin{thebibliography}{99}

\bibitem{Ab-H}
T. Abdellaoui and H. Heinich: ``Sur la distance de deux lois dans le
cas vectoriel''. CRAS, Paris S\'erie I, Math., {\bf 319}, 397-400, 1994.


\bibitem{An-K}
R. D. Anderson and V.L. Klee, Jr.: ``Convex functions and upper
semicontinuous collections''. Duke Math. Journal, {\bf 19}, 349-357, 1952.

\bibitem{App-1}
P. Appell: ``M\'emoire sur d\'eblais et les remblais des syst\`emes
continus ou discontinus''. M\'emoires pr\'esent\'ees par divers
savants \`a l'Acad\'emie des Sciences de l'Institut de France. Paris,
I. N. {\bf 29}, 1-208, 1887.

\bibitem{App-2}
P. Appell: ``Le probl\`eme g\'eom\'etrique des d\'eblais et des
remblais''. M\'emorial des Sciences Math\'ematiques, fasc. {\bf
XXVII}, Paris, 1928.


\bibitem{B-F}
P. J. Bickel and D. A. Freedman:
``Some asymptotic theory for the bootstrap''. The Annals of
Statistics, Vol. {\bf 9}, No. 6, 1196-1217, 1981.


\bibitem{BRE}
Y. Brenier: 
``Polar factorization and monotone rearrangement of vector valued
functions''. Comm. pure Appl. Math, {\bf 44}, 375-417, 1991.

\bibitem{Caf}
L. A. Caffarelli: ``The regularity of mappings with a convex
potential''. Jour. Amer. Math. Soc., {\bf 5}, 99-104, 1992.




\bibitem{D-M}
B. Dacarogna and J. Moser:
``On a partial differential equation involving the Jacobian
determinant''. Ann. Inst. Henri Poincar\'e, Analyse non-lin\'eaire,
{\bf 7}, 1-26, 1990.

\bibitem{Del-M}
C. Dellacherie and P. A. Meyer: {\sl Probabilit\'es et Potentiel,
  Ch. I \`a IV}. Paris, Hermann, 1975.

\bibitem{Du-S}
N. Dunford and J. T. Schwartz: {\sl Linear Operators} {\bf 2},
Interscience 1963.

\bibitem{Fer1}
X. Fernique:
``Extension du th\'eor\`eme de Cameron-Martin aux translations
al\'eatoires'', Comptes Rendus Math\'ematiques, Vol. {\bf 335}, Issue 1,
 65-68, 2002.   

\bibitem{Fer2}
X. Fernique: ``Comparaison aux mesures gaussiennes, espaces
autoreproduisants. Une application des propri\'et\'es
isop\'erim\'etriques''. Preprint.


\bibitem{F-P}
D. Feyel and A. de La Pradelle: ``Capacit\'es gaussiennes''. Annales
de l'Institut Fourier, t.41, f.1, 49-76, 1991.



\bibitem{F-U1}
D. Feyel and A. S. \"Ust\"unel: ``The notion of convexity and concavity
on Wiener space''. Journal of Functional Analysis, {\bf 176}, 400-428,
2000.

\bibitem{F-U2}
D. Feyel and A. S. \"Ust\"unel: ``Transport of measures on Wiener
space and the Girsanov theorem''. Comptes Rendus Math\'ematiques,
Vol. {\bf 334}, Issue 1, 1025-1028,  2002.

\bibitem{G-Mc}
W. Gangbo and R. J. McCann: ``The geometry of optimal
transportation''. Acta Mathematica, {\bf 177}, 113-161, 1996.

\bibitem{I-N}
K. Ito and M. Nisio: ``On the convergence of sums of independent
Banach space valued random variables''. Osaka Journ. Math. {\bf 5},
35-48, 1968.

\bibitem{Kan}
L. V. Kantorovitch: ``On the transfer of
masses''. Dokl. Acad. Nauk. SSSR {\bf 37}, 227-229, 1942.




\bibitem{Mar}
K. Marton: ``Bounding $\bar{d}$-distance by informational divergence:
a method to prove measure concentration''. Annals of Probability, {\bf
  24}, no.2, 857-866, 1996.


\bibitem{Mc1}
R. J. McCann: ``Existence and uniqueness of monotone measure-preserving
maps''. Duke Math. Jour., {\bf 80}, 309-323, 1995.

\bibitem{Mc2}
R. J. McCann: ``A convexity principle for interacting
gases''. Advances in Mathematics, {\bf 128}, 153-179, 1997.

\bibitem{Monge}
G. Monge: ``M\'emoire sur la th\'eorie des d\'eblais et des
remblais''. Histoire de l'Acad\'emie Royale des Sciences, Paris, 1781.

\bibitem{ROC}
T. Rockafellar: {\sl Convex Analysis}. Princeton University Press,
Princeton, 1972.

\bibitem{Sud}
V. N. Sudakov: ``Geometric problems in the theory of infinite
dimensional probability distributions''. Proc. Steklov Inst. Math.,
{\bf 141}, 1-178, 1979.





\bibitem{Tal}
M. Talagrand: ``Transportation cost for Gaussian and other product
measures''. Geom. Funct. Anal., {\bf 6}, 587-600, 1996.



\bibitem{THOMAS}
E. Thomas:
``The Lebesgue-Nikodym theorem for vector valued Radon
measures''. Memoirs of A.M.S., {\bf  139}, 1974.

\bibitem{ASU-0}
A. S. \"Ust\"unel: ``Representation of distributions on Wiener space
and Stochastic Calculus of Variations''. Journal of Functional
Analysis, {\bf 70}, 126-139, 1987.


\bibitem{ASU}
A. S. \"Ust\"unel:
{\sl Introduction to Analysis on Wiener Space}.
Lecture Notes in Math. Vol. {\bf 1610}. Springer, 1995.

\bibitem{BOOK}
A. S. \"Ust\"unel and M. Zakai:
{\sl Transformation of Measure on Wiener Space}.
Springer Monographs in Mathematics. Springer Verlag, 1999.
\end{thebibliography}
\end{document}